\documentclass[10pt,a4paper]{amsart}
\usepackage{amsmath}
\usepackage{amsfonts}
\usepackage{amsthm}
\usepackage{mathtools}
\usepackage{amscd}
\usepackage[utf8]{inputenc}
\usepackage{t1enc}
\usepackage[mathscr]{eucal}
\usepackage{indentfirst}
\usepackage{graphicx}
\usepackage{graphics}
\usepackage{pict2e}
\usepackage{epic}
\usepackage{float}
\usepackage{MnSymbol}
\usepackage{multirow}
\usepackage{shuffle} 
\usepackage[table,xcdraw]{xcolor} 
\usepackage{hyperref}
\usepackage{tcolorbox}
\numberwithin{equation}{section}
\usepackage[margin=2.9cm]{geometry}
\usepackage{epstopdf}
\usepackage{listings}
\allowdisplaybreaks

\newmuskip\pFqmuskip
\newcommand*\pFq[6][8]{%
	\begingroup 
	\pFqmuskip=#1mu\relax
	\mathchardef\normalcomma=\mathcode`,
	\mathcode`\,=\string"8000
	\begingroup\lccode`\~=`\,
	\lowercase{\endgroup\let~}\pFqcomma
	{}_{#2}F_{#3}{\left[\genfrac..{0pt}{}{#4}{#5};#6\right]}%
	\endgroup
}
\newcommand{\pFqcomma}{{\normalcomma}\mskip\pFqmuskip}

\theoremstyle{plain}
\newtheorem{theorem}{Theorem}[section]
\newtheorem{lemma}[theorem]{Lemma}

\newtheorem{proposition}[theorem]{Proposition}

\theoremstyle{definition}
\newtheorem{definition}[theorem]{Definition}

\newtheorem{remark}[theorem]{Remark}
\newtheorem{problem}[theorem]{Problem}

\newcommand{\CMZV}[2]{\textsf{CMZV}^{#1}_{#2}}

\DeclareMathOperator{\Li}{Li}

\newtheorem{romexample}{Example}

\usepackage{accsupp}
\DeclareRobustCommand\squelch[1]{%
	\BeginAccSupp{method=plain,ActualText={}}#1\EndAccSupp{}}

\newcounter{inprcnt}

\newcommand{\outpr}{\textcolor{gray!80!white}{\squelch{\smaller \texttt{Out[\arabic{inprcnt}]= }}}}

\usepackage{environ}

\makeatletter%
\NewEnviron{mathout}{%
	\small%
	\setlength\@mathmargin{0pt}%
	\setlength{\abovedisplayskip}{\bigskipamount}%
	\setlength{\belowdisplayskip}{\bigskipamount}%
	\setlength{\abovedisplayshortskip}{\bigskipamount}%
	\setlength{\belowdisplayshortskip}{\bigskipamount}%
	\begin{flalign}%
		\BODY
	\end{flalign}%
}
\makeatother%

\begin{document}
	
	\title{Wilf-Zeilberger seeds and non-trivial hypergeometric identities}
	
	\author[Kam Cheong Au]{Kam Cheong Au}
	
	\address{Rheinische Friedrich-Wilhelms-Universität Bonn \\ Mathematical Institute \\ 53115 Bonn, Germany} 
	
	\email{s6kmauuu@uni-bonn.de}
	

	\begin{abstract}Through a systematic approach on generating Wilf-Zeilberger-pairs, we prove some hypergeometric identities conjectures due to Z.W. Sun, J. Guillera and Y. Zhao etc., including two Ramanujan-$1/\pi^4$, one $1/\pi^3$ formulas as well as a remarkable series for $\zeta(5)$. 
	\end{abstract}
	
	\maketitle
	
	\section{Introduction}
	Ramanujan's $1/\pi$-formulas are infinite series of form
	$$\sum_{n\geq 0} \frac{(1/2)_n (1/2-s)_n (1/2+s)_n}{(1)_n} (an+b) z^n = \frac{c}{\pi}, \qquad a,b,c,z \in \overline{\mathbb{Q}}, \quad s\in \{0,1/3,1/4,1/6\}.$$
	They are now well-understood: they arise from the hypergeometric nature of modular forms: $a,b,c,z$ being algebraic comes from CM theory of modular and quasi-modular forms (\cite{borwein1987pi}, \cite{shimura1971introduction}, \cite{cohen2021rational}). J. Guillera brings the subject a new wind by discovering and proving generalizations of such formulas, for example (\cite{guillera2011new}), 
	\begin{equation}\tag{*}\label{pi-2_eq1}\sum_{n\geq 0} \frac{(1/2)_n^5}{(1)_n^5} \frac{(-1)^n}{2^{10n}} (820n^2+180n+13) = \frac{128}{\pi^2}.\end{equation}
	Such $1/\pi^k$-formulas have been more elusive than their $1/\pi$ counterparts. We now know a few examples of Ramanujan-type $1/\pi^{2}$-formulas (see \cite{guillera2016bilateral}). Some connections of them to Calabi-Yau differential equations and Hilbert modular forms are known (\cite{zudilin2011arithmetic}, \cite{zudilin2007quadratic}, \cite{dembele2022special}), indicating a deep underlying modular connection. \par
	However, up to now, the only known complete proof of equation (\ref{pi-2_eq1}) (and several others) seems to come from Wilf-Zeilberger-pair (WZ pair), a technique combinatorial in nature. \\[0.01in]
	
	A goal of this article is, after introducing the concept of Wilf-Zeilberger seed, to present a general method of manufacturing suitable WZ-pairs for proving a purported hypergeometric identity. Using this method, we establish following empirically-discovered $1/\pi^3, 1/\pi^4$ formulas\footnote{they are all the currently known $1/\pi^3, 1/\pi^4$ formulas, according to \cite{cohen2021rational}.}: 
$$\begin{aligned}\sum_{n\geq 0} \frac{1}{2^{6n}} \frac{(1/2)_n^7}{(1)_n^7}  (168n^3+76n^2+14n+1) &= \frac{32}{\pi^3}, \\
\sum _{n\geq 0} 2^{-12n} \frac{(\frac{1}{2})_n^7 (\frac{1}{4})_n(\frac{3}{4})_n}{(1)_n^9}(43680 n^4+20632 n^3+4340 n^2+466 n+21)&=\frac{2048}{\pi ^4}, \\
\sum _{n\geq 0} (-\frac{27}{256})^n \frac{(\frac{1}{2})_n^5 (\frac{1}{4})_n(\frac{3}{4})_n(\frac{1}{3})_n(\frac{2}{3})_n}{(1)_n^9}(4528 n^4+3180 n^3+972 n^2+147 n+9)&=\frac{768}{\pi ^4},\end{aligned}$$
	
	due to B. Gourevich, J. Cullen and Y. Zhao respectively (Examples \ref{ex_pi-3}, \ref{ex_pi-4_1}, \ref{ex_pi-4_2}). 
	
	This paper can be seen as a sequel to the author's previous article \cite{au2022multiple}. Actually, Ramanujan $1/\pi^c$-series does not occupy a special place in the method, it can be applied to many series of the form \begin{equation}\label{to_be_proved_identity_0}\sum_{n\geq 0} z^n \frac{(a_1)_n\cdots (a_m)_n}{(b_1)_n \cdots (b_m)_n} R(n) \in \frac{\text{CMZVs of some level}}{\pi^c}, \end{equation}
	here CMZVs stand for colored multiple zeta values (\cite{zhao2016multiple}), and $$z\in \mathbb{Q},\quad a_i,b_i\in \mathbb{Q}\cap (0,1], \quad R(n)\in \mathbb{Q}(n), \quad c\in \mathbb{N}.$$
	Some other formulas that we will establish include
	 $$\sum_{n\geq 1}\left(-\frac{2^{10}}{5^5}\right)^n \frac{(1)_n^9}{(\frac{1}{2})_n^5(\frac{1}{5})_n(\frac{2}{5})_n(\frac{3}{5})_n(\frac{4}{5})_n} \frac{5532n^4-5600n^3+2275n^2-425n+30}{n^9} =-380928\zeta(5),$$
	 due to Y. Zhao (Example \ref{ex_zeta5}) and 
	 $$\sum_{n\geq 0} (\frac{4}{27})^n \frac{(1/2)_n^7}{(1/6)_n (5/6)_n (1)^5} \frac{92n^3+54n^2+12n+1}{6n+1} = \frac{12}{\pi^2},$$
	 due to Z.W. Sun (Example \ref{12pi-2_ex}). We also include many conjectures of Z.W. Sun (\cite{sun2023new}, \cite{sun2021book}, \cite{sun2022conjectures}) to demonstrate versatility of the approach. \par 
	 
	 Nonetheless, this approach of proving infinite sum via WZ-pairs is still inchoate and more developments are necessary if we want to use it to prove more non-trivial identities like
	$$\sum_{n\geq 0} (\frac{1}{7^{4}})^n \frac{(1/2)_n (1/8)_n (3/8)_n (5/8)_n (7/8)_n}{(1)_n^5} (1920n^2+304n+15) \stackrel{?}{=} \frac{56\sqrt{7}}{\pi^2},$$
found experimentally by J. Guillera in \cite{guillera2003new}.
	~\\[0.03in]

We structure the paper as follows: section 2 provides, as two motivating examples, complete proof of two $1/\pi^4$-formulas; section 3 introduces the key concept of WZ-seed, a helpful concept on producing WZ-pairs; section 4 and 5 consists of examples using this method of producing WZ-pair.

	\section{The central proposition and two motivating examples}
	The following proposition is an important bridge between WZ-pairs and infinite series. 
	\begin{proposition}\label{WZ_prop}
	Suppose $F,G: \mathbb{N}^2\to \mathbb{C}$ are two functions such that
	$$F(n+1,k) - F(n,k) = G(n,k+1)-G(n,k)$$
	and \begin{itemize}
		\item $\sum_{k\geq 0} F(0,k)$ converges,
		\item $\sum_{n\geq 0} G(n,0)$ converges,
		\item $\lim_{k\to\infty} G(n,k) = g(n)$ exists for each $n\in \mathbb{N}$ and $\sum_{n\geq 0} g(n)$ converges,
		\end{itemize}
	then $\lim_{n\to \infty} \sum_{k\geq 0} F(n,k)$ exists and is finite, also
	$$\sum_{k\geq 0} F(0,k) + \sum_{n\geq 0} g(n) = \sum_{n\geq 0} G(n,0) + \lim_{n\to \infty} \sum_{k\geq 0} F(n,k).$$
	\end{proposition}
\begin{proof}
	Apply $\sum_{k=0}^{K-1}\sum_{n=0}^{N-1}$ on both sides of $F(n+1,k) - F(n,k) = G(n,k+1)-G(n,k)$, one gets
	$$\sum_{k=0}^{K-1} F(N,k) - \sum_{k=0}^{K-1} F(0,k) = \sum_{n=0}^{N-1} G(n,K) - \sum_{n=0}^{N-1} G(n,0),$$
	by assumption, we can let $K\to\infty$ while fixing $N$, arriving at
	$$\sum_{k\geq 0} F(N,k) - \sum_{k\geq 0} F(0,k) = \sum_{n=0}^{N-1} g(n) - \sum_{n=0}^{N-1} G(n,0),$$
	letting $N\to \infty$ gives the assertion. 
	\end{proof}

We first look at two examples of Ramanujan-type $1/\pi^4$ formulas.

	\begin{romexample}\label{ex_pi-4_1}
		Let us establish the following Ramanujan-$1/\pi^4$ series due to J. Cullen in 2010 (\cite{guillera2016bilateral}):
		$$\sum _{n\geq 0} 2^{-12n} \frac{(\frac{1}{2})_n^7 (\frac{1}{4})_n(\frac{3}{4})_n}{(1)_n^9}(43680 n^4+20632 n^3+4340 n^2+466 n+21)=\frac{2048}{\pi ^4}.$$
		This is an easy application of above proposition once we somehow write down a suitable $F(n,k)$. Let $$F(n,k) = \frac{\left(2 k+3 n+\frac{1}{2}\right) \Gamma (1-2 n)^3 \Gamma \left(n+\frac{1}{2}\right)^3 \Gamma \left(k+n+\frac{1}{2}\right)^5 \Gamma (k+4 n)}{\Gamma \left(\frac{1}{2}-n\right)^7 \Gamma (2 n)^2 \Gamma \left(k-n+\frac{3}{2}\right) \Gamma (k+2 n+1)^5},$$
		one checks it has a WZ-mate\footnote{readers interested in full expression should perform Gosper's algorithm themselves.}
		$$G(n,k) = F(n,k)\frac{(-2 k+2 n-1) (68728320 k^2 n^8 + 41932800 k n^9 + 11182080 n^{10} + \text{59 more monomials} )}{131072 n^5 (k+2 n+1)^5 (4 k+6 n+1)}.$$
		It is easy to see\footnote{we will elaborate more about this point in next section} $\lim_{n\to\infty} \sum_{k\geq 0} F(n,k) = 0$, and $\lim_{k\to\infty} G(n,k) := g(k)$ exists and is non-zero. For a generic $a\in \mathbb{C}$, we have $$\sum_{k\geq 0} F(a,k) + \sum_{n\geq 0} g(n+a) = \sum_{n\geq 0} G(n+a,0),$$
		which becomes, after dividing both sides by $\frac{4 a^2 \Gamma (1-2 a)^3 \Gamma (a+1/2)^8 \Gamma (4 a+1)}{\Gamma(1/2-a)^8 \Gamma (2 a+1)^7}$,
		\begin{multline*}\sum_{k\geq 0}\frac{(6 a+4 k+1) (a+1/2)_k^5 (4 a+1)_k}{(1+2k-2a) (4 a+k) (1/2-a)_k (2 a+1)_k^5} \\ - \sum_{n\geq 0} \frac{(-1)^n (820 a^2+1640 a n+180 a+820 n^2+180 n+13) \Gamma (1/2-a) \Gamma (2 a+1)^5 (a+\frac{1}{2})_n^{10}}{4096 a^5 \Gamma (a+1/2)^5 \Gamma(4 a+1)(2 a+1)_{2 n}^5} = \\ -\sum_{n\geq 0} \frac{\splitfrac{(21 + 466 a + 4340 a^2 + 20632 a^3 + 43680 a^4 + 
			466 n + 8680 a n + 61896 a^2 n + 174720 a^3 n +}{  4340 n^2 + 
			61896 a n^2 + 262080 a^2 n^2 + 20632 n^3 + 174720 a n^3 + 
			43680 n^4) (1/2 + a)_n^{16} (1 + 4 a)_{ 
			4 n}}}{65536 a^5 (1 + 2 a)_{2 n}^{10}}.
			\end{multline*}
		By comparing coefficient of $a^{-5}$ in Laurent expansion at $a=0$  , to which first term of LHS does not contribute, we have
		$$\frac{1}{4096\pi^2} \sum_{n\geq 0} (-1)^n \frac{(1/2)_n^{10}}{(1)_{2n}^5} (13+180n+820n^2)  = \sum_{n\geq 0} \frac{(1/2)_n^{16} (1)_{4n}}{65536(1)_{2n}^{10}} (43680 n^4+20632 n^3+4340 n^2+466 n+21).$$
		RHS is exactly the infinite series of our concern; for LHS, it is already shown by J. Guillera (\cite{guillera2011new}), also using WZ-pair, that $$\sum_{n\geq 0} \frac{(1/2)_n^5}{(1)_n^5} \frac{(-1)^n}{2^{10n}} (820n^2+180n+13) = \frac{128}{\pi^2}.$$ This completes the proof. 
		
		\end{romexample}

	\begin{romexample}\label{ex_pi-4_2}The following series, conjectured by Y. Zhao\footnote{\url{https://mathoverflow.net/questions/281009}}, 
	$$\sum _{n\geq 0} (-\frac{27}{256})^n \frac{(\frac{1}{2})_n^5 (\frac{1}{4})_n(\frac{3}{4})_n(\frac{1}{3})_n(\frac{2}{3})_n}{(1)_n^9}(4528 n^4+3180 n^3+972 n^2+147 n+9)=\frac{768}{\pi ^4},$$
	can be proved \textit{mutatis mutandis} as above example with another $F(n,k)$: $$F(n,k) = \frac{(-1)^n (1/2+2k+3n) \Gamma (1-n)^3 \Gamma \left(n+\frac{1}{2}\right)^3 \Gamma \left(k+n+\frac{1}{2}\right)^4 \Gamma \left(k+2 n+\frac{1}{2}\right) \Gamma (k+3 n)}{\Gamma \left(k+\frac{3}{2}\right) \Gamma \left(\frac{1}{2}-n\right)^3 \Gamma (n) \Gamma (2 n) \Gamma (k+n+1) \Gamma (k+2 n+1)^4}.$$ 
	then, after finding $G(n,k), g(n)$, $\sum_{k\geq 0} F(a,k) + \sum_{n\geq 0} g(n+a) = \sum_{n\geq 0} G(n+a,0)$ becomes
	\begin{multline*}\sum_{k\geq 0} \frac{(6 a+4 k+1) (1/2+a)_k^4 (1/2+2a)_k (1+3a)_k}{(2 k+1) (3 a+k) (1+a)_k (1+2a)_k^4} - \\ \sum_{n\geq 0} \frac{(-1)^n (1 + 8 a + 20 a^2 + 8 n + 40 a n + 20 n^2) \sqrt{\pi} \Gamma(1+a)\Gamma(1+2a)^4 (1/2+a)_n^6}{16a^5 \Gamma(1/2+a)^4 \Gamma(1/2+2a) \Gamma(1+3a) (1+a)_n^4 (1+2a)_{2n}} \\ = - \sum_{n\geq 0} \frac{\splitfrac{(-1)^n (9 + 147 a + 972 a^2 + 3180 a^3 + 4528 a^4 + 147 n + 
			1944 a n + 9540 a^2 n + 18112 a^3 n + 972 n^2 + 9540 a n^2}{ + 
			27168 a^2 n^2 + 3180 n^3 + 18112 a n^3 + 4528 n^4) (1/2+a)_n^{10} (1/2+2a)_{2n} (1+3a)_{3n}}}{1536a^5 (1+a)_n^5 (1+2a)_{2n}^5}. \end{multline*}
Comparing coefficient of $a^{-5}$ gives
		$$\sum_{n\geq 0} \frac{(-1)^n (1+8n+20n^2) (1/2)_n^6}{16\pi^2 (1)_n^4 (1)_{2n}} = \sum_{n\geq 0} \frac{(-1)^n (9+147n+972n^2+3180n^3+4528n^4) (1/2)_n^{10} (1/2)_{2n} (1)_{3n}}{1536 (1)_n^5 (1)_{2n}^5},$$
		RHS is exactly the infinite series of our concern, while LHS is a known result in \cite{guillera2011new} $$\sum_{n\geq 0} \frac{(1/2)_n^5}{(1)_n^5} \frac{(-1)^n}{2^{2n}} (20n^2+8n+1) = \frac{8}{\pi^2}.$$ This completes the proof. 
	\end{romexample}

In next section, we will explain how did we come up with $F(n,k)$ in above two examples.

	\section{Wilf-Zeilberger seeds}
	\subsection{A source of WZ-pairs}
	Following definition is paramount in our methodology of finding suitable WZ-pairs:
	\begin{definition}
	Let $f(a_1,\cdots,a_m,k)$ be a hypergeometric term in $a_i$ and $k$, it is called a \textit{WZ-seed} in variables $a_i$ if for all $A_i \in \mathbb{Z}$ and $K\in \mathbb{Z}$,
	$$F(n,k) = f(a_1+A_1 n,a_2+A_2 n,\cdots,a_m+A_m n,k+Kn)$$
	has a hypergeometric WZ-mate $G(n,k)$. 
	\end{definition}

Following terminology of \cite{gessel1995finding}, we call $a_i$ the \textit{accessory parameters} of $F(n,k)$. 

The following empirical way to produce (potential) WZ-seeds are known at least back to Zeilberger, Gessel \cite{gessel1995finding} even used it systematically to generate many terminating hypergeometric identities. 
	\begin{problem}[Hypergeometric summation formula induces WZ-candidate]\label{summation_induces_WZ}
	Let $f(a_1,\cdots,a_m;k)$ be a proper hypergeometric term in $a_i,k$. If $$\sum_{k\geq 0} f(a_1,\cdots,a_m,k)$$ exists and is independent of $a_i$, then is it true that, $f(a_1,\cdots,a_m;k)$ is a WZ-seed?
\end{problem}

The answer seems to be affirmative in all examples we have computed. Gessel (\cite{gessel1995finding}) already used this empirical observation to generate a plethora of finite hypergeometric identities. We give an example of above observation, consider the Gauss $_2F_1$ formula:
$$\pFq{2}{1}{a \quad b}{c}{1} = \sum_{k\geq 0} \frac{(a)_k(b)_k}{(1)_k(c)_k} = \frac{\Gamma(c-a-b)\Gamma(c)}{\Gamma(c-a)\Gamma(c-b)},$$
divide the gamma factor of RHS into LHS summand, we arrive at:
$$\sum_{k\geq 0} \frac{\Gamma (c-a) \Gamma (a+k) \Gamma (c-b) \Gamma (b+k)}{\Gamma (a) \Gamma (b) \Gamma (k+1) \Gamma (c+k) \Gamma (-a-b+c)} = 1.$$
If we denote the summand by $\textsf{Gauss2F1}(a,b,c,k)$, then the above observation would predict $$F(n,k) = \textsf{Gauss2F1}(a - n,b, c+n, k+n)$$
has a WZ-mate $G$ that is also proper hypergeometric, indeed, Gosper's algorithm returns
\begin{multline*}G(n,k) =F(n,k) \times (-a^2 b+a^2 c+a^2 n+2 a b c+a b k+5 a b n+a b-2 a c^2-2 a c k-8 a c n-a c-a k^2-4 a k n-7 a n^2-a n\\ -b c^2-2 b c k-6 b c n-b c-b k^2-5 b k n-8 b n^2-2 b n+c^3+2 c^2 k+7 c^2 n+c^2+c k^2+8 c k n+c k+15 c n^2+4 c n+2 k^2 n+k^2\\ +8 k n^2+3 k n+10 n^3+4 n^2) \div((c+k+2 n) (a+b-c-2 n-1) (a+b-c-2 n)), \end{multline*}
it satisfies $F(n+1,k)-F(n,k) = G(n,k+1)-G(n,k)$. \par


~\\[0.02in]
To check a given hypergeometric term is indeed a WZ-seed, the definition is not practical because it requires existence of WZ-mate for infinitely many choices of $A_1,\cdots,A_m \in \mathbb{Z}$ and $K\in \mathbb{Z}$. We now give a more conceptual description of WZ-seed. 
\begin{theorem}[Criterion of being WZ-seed]\label{WZ_crit}
Let $f(a_1,\cdots,a_m,k)$ be a hypergeometric term in $a_i$ and $k$. The following are equivalent:
\begin{enumerate}
	\item $f(a_1,\cdots,a_m,k)$ is a WZ-seed with accessory parameters $a_i$;
	\item for $i=1,\cdots,m$, there exists a rational function $r_i$ such that $\Delta_{a_i}(f) = \Delta_k (r_i f)$.\footnote{Here $\Delta_x$ means forward difference with respect to variable in the subscript: $\Delta_x f(x) = f(x+1)-f(x)$}
\end{enumerate}
\end{theorem}
\begin{proof}
The implication $(1)\implies (2)$ is evident: $(1)$ implies there exists a rational function $r_i$ such that
$$f(a_1,\cdots,a_i+n + 1,\cdots,a_m,k) - f(a_1,\cdots,a_i+n ,\cdots,a_m,k) = \Delta_n (r_i f).$$
Because the variable $a_i$ always appears with $n$ in form of $a_i+n$, we have $\Delta_n (r_i f) = \Delta_{a_i} (r_i f)$. Setting $n=0$ gives
$$f(a_1,\cdots,a_i+ 1,\cdots,a_m,k) - f(a_1,\cdots,a_i ,\cdots,a_m,k) = \Delta_{a_i} (r_i f),$$
showing $(2)$. Now we show $(2)\implies (1)$, let $A_i \in \mathbb{Z}, K\in \mathbb{Z}$ be given. Denote $\textbf{S}_{a_i}$ to be shifting operator that changes $a_i$ to $a_i+1$. For integer $l$, $\textbf{S}_{a_i}^l$ acts by changing $a_i$ to $a_i+l$. $(2)$ means there exists $r_i$ such that
$$(\textbf{S}_{a_i}-1)f = \Delta_k (r_i f), \qquad 1\leq i\leq m.$$
In (commutative) ring $R = \mathbb{Q}[\textbf{S}_{a_1},\cdots,\textbf{S}_{a_m},\textbf{S}_k,\textbf{S}_{a_1}^{-1},\cdots,\textbf{S}_{a_m}^{-1},\textbf{S}_k^{-1}]$, ideal $(\textbf{S}_{a_1}-1,\cdots,\textbf{S}_{a_m}-1,\textbf{S}_k-1)$ is maximal, it contains any rational function in $R$ that vanishes at point $(\textbf{S}_{a_1},\cdots,\textbf{S}_{a_m},\textbf{S}_k)=(1,\cdots,1)$. $\left(\prod_{i=1}^m \textbf{S}_{a_i}^{A_i}\right)\textbf{S}_k^K - 1$ is in this ideal, so there exists $\textbf{g}_i \in R$ and $\textbf{h}\in R$ such that
$$\left(\prod_{i=1}^m \textbf{S}_{a_i}^{A_i}\right)\textbf{S}_k^K - 1 = \sum_{i=1}^m \textbf{g}_i (\textbf{S}_{a_i}-1) + \textbf{h} (\textbf{S}_k-1),$$
Applying this to $f(a_1,\cdots,a_m,k)$ gives
$$f(a_1+A_1,\cdots,a_m+A_m,k+K)-f(a_1,\cdots,a_m,k) = \sum_{i=1}^m \textbf{g}_i \Delta_k(r_i f) + \textbf{h}(\textbf{S}_k-1)f= \Delta_k \left(\sum_{i=1}^m \textbf{g}_i r_i f + \textbf{h}f \right)$$
Replace $a_i$ by $a_i + A_i n$, $k$ by $k+Kn$, LHS becomes
$$\Delta_n (f(a_1+A_1n,\cdots,a_m+A_m n,k+Kn))$$
so we see $F(n,k) = f(a_1+A_1n,\cdots,a_m+A_m n,k+Kn)$ has a WZ-mate given by $$(\textbf{g}_i r_i +\textbf{h}) f(a_1+A_1n,\cdots,a_m+A_m n, k+Kn)$$
this completes the proof of implication $(2)\implies (1)$. 
\end{proof}

Using above theorem, one can demonstrate $f(a,b,c,k) = \textsf{Gauss2F1}(a,b,c,k)$ is indeed a WZ-seed. Invoking Gosper's algorithm on variables $a,b,c$, we find
$$\Delta_a(f) = \Delta_k(r_1 f), \qquad \Delta_b(f) = \Delta_k(r_2 f), \qquad \Delta_c(f) = \Delta_k(r_3 f),$$
with \begin{equation}\label{Gauss2F1abc}r_1 = \frac{k (c+k-1)}{a (a-c+1)} ,\qquad 
r_2 = \frac{k (c+k-1)}{b (b-c+1)}, \qquad
r_3 = -\frac{k}{a+b-c}.\end{equation}

\subsection{Computational implementation}
Implementations of Gosper's algorithm are abundant, the one we use in this article comes from P. Paule, M. Schorn, and A. Riese (\cite{paule1995mathematica}), using programming language of Mathematica. Readers are referred to appendix for details of implementation. \par

First we load the Mathematica package, which can be obtained from \url{https://www3.risc.jku.at/research/combinat/software/ergosum} as a part of larger RISCErgoSum package.
\begin{lstlisting}[]
<@\inpr@> Needs["RISC`fastZeil`"];
\end{lstlisting}

The package contains a command \texttt{Gosper} which enables us to find WZ-mate:
\begin{lstlisting}[]
<@\inpr@> findWZmate[f_,n_,k_] := Module[{gosper},
 gosper=Gosper[(FunctionExpand[(f/.{n->n+1})/f]-1)*f,k];
 If[gosper!={},gosper[[1,2,2]]/f,$Failed]];
\end{lstlisting}

Given a $F(n,k)$, \texttt{findWZmate} finds rational function $R(n,k)$, if exists, such that for $G(n,k) = F(n,k)R(n,k)$, we have $\Delta_n F = \Delta_k G$. For example, using our example of seed $\textsf{Gauss2F1}$

\begin{lstlisting}[]
<@\inpr@> Gauss2F1[a_,b_,c_,k_]:=Gamma[-a+c]Gamma[-b+c]Gamma[a+k]Gamma[b+k]/
(Gamma[a]Gamma[b]Gamma[-a-b+c]Gamma[1+k]Gamma[c+k]);
\end{lstlisting}

following command finds WZ-mate of $\textsf{Gauss2F1}(a-n,b,c+n,k)$:
\begin{lstlisting}[]
<@\inpr@> findWZmate[Gauss2F1[a-n,b,c+n,k],n,k]
\end{lstlisting}
\begin{mathout}
\tag*{\outpr}\frac{k \left(-a^2-a k+2 a n+2 a-b c-b k-b n+b+c^2+c k+2 c n-c+2 k n+k-3 n-1\right)}{(a+k-n-1) (a+b-c-2 n-1) (a+b-c-2 n)} &&
\end{mathout}

and following finds WZ-mate $G(n,k)$ for the $F(n,k)$ in Example \ref{ex_pi-4_1}:
\begin{lstlisting}[]
<@\inpr@> findWZmate[(2k+3n+1/2)Gamma[1-2n]^3*Gamma[n+1/2]^3*
Gamma[k+n+1/2]^5*Gamma[k+4n]/(Gamma[1/2-n]^7*Gamma[2n]^2
Gamma[k-n+3/2]*Gamma[k+2n+1]^5),n,k]
\end{lstlisting}

The rational functions $r_1,r_2,r_3$ in equation (\ref{Gauss2F1abc}) are found using
\begin{lstlisting}[]
<@\inpr@> findWZmate[Gauss2F1[a,b,c,k],#,k] & /@ {a,b,c}
\end{lstlisting}
thereby confirming $\textsf{Gauss2F1}(a,b,c,k)$ is a WZ-seed by Theorem \ref{WZ_crit}. 

	\subsection{List of WZ-seeds}
	
We exhibit some WZ-seeds (one can verify this by Theorem \ref{WZ_crit}). The following list is by no mean exhaustive, many seeds which do not generate useful examples are not included. Most of this list are already known to Gessel \cite{gessel1995finding}.

\begin{enumerate}
\item  $\textsf{Gauss2F1}(a,b,c,k)$
$$ = \frac{\Gamma (c-a) \Gamma (a+k) \Gamma (c-b) \Gamma (b+k)}{\Gamma (a) \Gamma (b) \Gamma (k+1) \Gamma (c+k) \Gamma (-a-b+c)}.$$
Origin: Gauss $_2F_1$ summation formula. Used in: Example \ref{gauss_ex1}; Examples I, II, V of \cite{au2022multiple}.

\item  $\textsf{Dixon3F2}(a,b,c,k)$
$$ = \frac{\Gamma (a-b+1) \Gamma (a-c+1) \Gamma (2 a+k) \Gamma (b+k) \Gamma (c+k) \Gamma (2 a-b-c+1)}{\Gamma (a) \Gamma (b) \Gamma (c) \Gamma (k+1) \Gamma (a-b-c+1) \Gamma (2 a-b+k+1) \Gamma (2 a-c+k+1)}.$$
Origin: Dixon's $_3F_2$ summation formula. Used in: Examples III, IV of \cite{au2022multiple}.

\item $\textsf{Watson3F2}(a,b,c,k)$
$$ = \frac{2^{-2 a-2 b+2 c} \Gamma \left(-a+c+\frac{1}{2}\right) \Gamma (2 a+k) \Gamma \left(-b+c+\frac{1}{2}\right) \Gamma (2 b+k) \Gamma (c+k)}{\Gamma (a) \Gamma (b) \Gamma (k+1) \Gamma (2 c+k) \Gamma \left(-a-b+c+\frac{1}{2}\right) \Gamma \left(a+b+k+\frac{1}{2}\right)}.$$ 
Origin: Watson's $_3F_2$ summation formula. Used in: Example IX of \cite{au2022multiple}.

\item $\textsf{Dougall5F4}(a,b,c,d,k)$
$$ = \frac{(a+2 k) \Gamma (a+k) \Gamma (b+k) \Gamma (c+k) \Gamma (d+k) \Gamma (a-b-c+1) \Gamma (a-b-d+1) \Gamma (a-c-d+1)}{\Gamma (b) \Gamma (c) \Gamma (d) \Gamma (k+1) \Gamma (a-b+k+1) \Gamma (a-c+k+1) \Gamma (a-d+k+1) \Gamma (a-b-c-d+1)}.$$ 
Origin: very-well-poised $_5F_4$ summation formula. Used in: Examples \ref{pi-2_parameters} and \ref{ex_reg}. Examples VI, VII, VIII of \cite{au2022multiple}.

\item $\textsf{Dougall7F6}(a,b,c,d,e,k)$
$$ = \frac{\splitfrac{(-1)^{a+e} (a+2 k) \Gamma (a+k) \Gamma (b+k) \Gamma (c+k) \Gamma (d+k) \Gamma (e+k) \Gamma (a-b-c+1) \Gamma (a-b-d+1)}{ \Gamma (a-c-d+1) \Gamma (-a+b+c+e) \Gamma (-a+b+d+e) \Gamma (-a+c+d+e) \Gamma (2 a-b-c-d-e+k+1)}}{\splitfrac{\Gamma (b) \Gamma (c) \Gamma (d) \Gamma (e) \Gamma (k+1) \Gamma (-a+b+e) \Gamma (a-b+k+1) \Gamma (-a+c+e) \Gamma (a-c+k+1) \Gamma (-a+d+e) \Gamma (a-d+k+1)}{ \Gamma (a-e+k+1) \Gamma (a-b-c-d+1) \Gamma (2 a-b-c-d-e+1) \Gamma (-a+b+c+d+e+k)}}.$$ 
Origin: very-well-poised $2$-balanced terminating $_7F_6$ summation formula. Used in: Examples \ref{ex_pi-4_1}, \ref{ex_pi-4_2}, \ref{ex_pi-3}

\item $\textsf{Balanced3F2}(a,b,c,d,k)$
$$ = \frac{\Gamma (d-a) \Gamma (a+k) \Gamma (d-b) \Gamma (b+k) \Gamma (d-c) \Gamma (c+k) (-1)^{a+b+c-d}}{\Gamma (a) \Gamma (b) \Gamma (c) \Gamma (k+1) \Gamma (d+k) \Gamma (-a-b+d) \Gamma (-a-c+d) \Gamma (-b-c+d) \Gamma (a+b+c-d+k+1)}.$$ 
Origin: $1$-balanced terminating $_3F_2$ summation formula.

\item $\textsf{Seed1}(a,c,k)$
$$ = \frac{2^{a-2 c-2 k} \Gamma (-a+2 c-1) \Gamma (a+2 k)}{\Gamma (a) \Gamma (k+1) \Gamma \left(-a+c-\frac{1}{2}\right) \Gamma (c+k)}.$$ 
Origin: put $a \to a/2, b\to a/2+1/2$ in \textsf{Gauss2F1}. Used in Examples \ref{-2beta4_ex} and \ref{ex_nonvanishing_1}. 

\item $\textsf{Seed2}(a,c,d,k)$
$$ = \frac{2^{-2 c-2 k} (-1)^{a+c-d} \Gamma (-a+2 d-1) \Gamma (a+2 k) \Gamma (d-c) \Gamma (c+k)}{\Gamma (a) \Gamma (c) \Gamma (k+1) \Gamma \left(-a+d-\frac{1}{2}\right) \Gamma (d+k) \Gamma (-a-2 c+2 d-1) \Gamma \left(a+c-d+k+\frac{3}{2}\right)}.$$ 
Origin: put $a \to a/2, b\to a/2+1/2$ in \textsf{Balanced3F2}. Used in Example \ref{384pi-2_ex}.

\item $\textsf{Seed3}(a,b,d,k)$
$$ = \frac{2^{2 d} (a+2 k) \Gamma \left(a-b+\frac{1}{2}\right) \Gamma (a+k) \Gamma (b+2 k) \Gamma (d+k) \Gamma (2 a-b-2 d+1)}{\Gamma (b) \Gamma (d) \Gamma (k+1) \Gamma \left(a-b-d+\frac{1}{2}\right) \Gamma (2 a-b+2 k+1) \Gamma (a-d+k+1)}.$$ 
Origin: put $b \to b/2, c\to b/2+1/2$ in \textsf{Dougall5F4}. Used in Examples \ref{896-zeta3_ex}, \ref{sun_100_prob_1}. 

\item $\textsf{Seed4}(a,b,k)$
$$ = \frac{3^{b+1} (a+2 k) \Gamma (3 a-2 b) \Gamma (a+k) \Gamma (b+3 k)}{\Gamma (b) \Gamma (k+1) \Gamma (a-b) \Gamma (3 a-b+3 k+1)}.$$ 
Origin: put $b \to b/3, c\to b/3+1/3, d\to b/3+2/3$ in \textsf{Dougall5F4}. Used in Examples \ref{ex_CMZVtrans}, \ref{ex_nonvanishing_2}. 

\item $\textsf{Seed7}(a,b,d,e,k)$
$$ = \frac{\splitfrac{(-1)^{a+e} (a+2 k) \Gamma \left(a-b+\frac{1}{2}\right) \Gamma (a+k) \Gamma (b+2 k) \Gamma (d+k) \Gamma (e+k) \Gamma (2 a-b-2 d+1)}{\Gamma \left(-a+b+e+\frac{1}{2}\right) \Gamma (-2 a+b+2 d+2 e) \Gamma \left(2 a-b-d-e+k+\frac{1}{2}\right)}}{\splitfrac{\Gamma (b) \Gamma (d) \Gamma (e) \Gamma (k+1) \Gamma \left(a-b-d+\frac{1}{2}\right) \Gamma (-2 a+b+2 e) \Gamma (2 a-b+2 k+1) \Gamma (-a+d+e) \Gamma (a-d+k+1)}{\Gamma (a-e+k+1) \Gamma \left(2 a-b-d-e+\frac{1}{2}\right) \Gamma \left(-a+b+d+e+k+\frac{1}{2}\right)}}.$$ 
Origin: put $b \to b/2, c\to b/2+1/2$ in \textsf{Dougall7F6}. Used in Examples \ref{ex_pi4_1}, \ref{12pi-2_ex}, \ref{ex_zeta5}. 

\item $\textsf{Seed9}(a,b,d,k)$
$$ = \frac{\splitfrac{(a+2 k) (-1)^d 2^{4 a-2 b-2 d} \Gamma \left(a-b+\frac{1}{2}\right) \Gamma (a+k) \Gamma (b+2 k) \Gamma (d+2 k)}{ \Gamma (-2 a+2 b+d+1) \Gamma (-2 a+b+2 d+1) \Gamma (2 a-b-d+k)}}{\Gamma (b) \Gamma (d) \Gamma (k+1) \Gamma \left(-a+d+\frac{1}{2}\right) \Gamma (-2 a+b+d) \Gamma (2 a-b+2 k+1) \Gamma (2 a-d+2 k+1) \Gamma (-a+b+d+k+1)}.$$ 
Origin: put $b \to b/2, c\to b/2+1/2, d\to d/2, e\to d/2+1/2$ in \textsf{Dougall7F6}. 

\item $\textsf{Seed10}(a,d,k)$
$$ = \frac{(-1)^{a-d} 3^{-a-3 k} \Gamma (-a+3 d-2) \Gamma (a+3 k)}{\Gamma (a) \Gamma (k+1) \Gamma (-2 a+3 d-3) \Gamma (d+k) \Gamma (a-d+k+2)}.$$ 
Origin: put $a\to a/3, b\to a/3+1/3, c\to a/3+2/3$ in \textsf{Balanced3F2}. Used in Example \ref{sun_100_prob_2}. 
\end{enumerate}

Many of the above seeds have already been used tacitly by other researchers. In Chu and Zhang \cite{chu2014accelerating}, an entire paper is devoted to formulas derived using WZ-seed \textsf{Dougall5F4}; in fact, their methods are equivalent to \ref{WZ_prop} used on WZ-seed \textsf{Dougall5F4}, a connection already spotted by J. Guillera in \cite{guillera2018dougall}. Zhang also wrote on formulas from seeds \textsf{Watson3F2} in \cite{zhang2015common} and \textsf{Gauss2F1} in \cite{chu2011dougall}. \par 

WZ-pairs employed in many literatures (\cite{mohammed2005infinite}, \cite{mohammed2005infinite}, \cite{pilehrood2008generating}, \cite{pilehrood2010series}, \cite{pilehrood2008simultaneous}) actually have their origins in seeds \textsf{Gauss2F1}, \textsf{Dixon3F2} and \textsf{Dougall5F4}. Many other works that try to prove special instances of equation (\ref{to_be_proved_identity}) actually uses formula that are derivable by WZ-pairs. (e.g. \cite{wei2023two}, \cite{wei2022conjectural}, \cite{wei2023some}, \cite{wei2023some2}, \cite{hou2023taylor}) \par

We remark that many WZ-seeds above have $q$-analogues and hence many examples in this paper have $q$-analogues as well. This is discussed in the author's another article \cite{au2024wilfQ}. 

\subsection{Searching a suitable WZ-pair}
Suppose we want to prove a numerical identity of form
\begin{equation}\label{to_be_proved_identity}\sum_{n\geq 0} z^n \frac{(a_1)_n\cdots (a_m)_n}{(b_1)_n \cdots (b_m)_n} R(n) \in \frac{\text{CMZVs of some level}}{\pi^r},  \qquad z\in \mathbb{C},\quad a_i,b_i\in \mathbb{Q}\cap (0,1], \quad R(n)\in \mathbb{Q}(n).\end{equation}
Our basic principal would be trying to search $F(n,k)$ such that for its WZ-mate $G(n,k)$, above series will be the same as $\sum_{n\geq 0} G(n,0)$. In all cases computed, the rational function $R(n)$ is irrelevant for the search, because it will always pop up automatically when running Gosper's algorithm. More crucial are the exponent $z$ and the pochhammer part $\frac{(a_1)_n\cdots (a_m)_n}{(b_1)_n \cdots (b_m)_n}$. \par

It is elementary to observe that exponent and pochhammer part of $G(n,0)$ can already be determined from expression of $F(n,k)$: indeed, since $G(n,k)/F(n,k)$ is rational, they are the same as exponent and pochhammer part of $F(n,0)$. 

For example, consider the seed labeled \textsf{Seed3} above, $$\textsf{Seed3}(a,b,d,k)= \frac{2^{2 d} (a+2 k) \Gamma \left(a-b+\frac{1}{2}\right) \Gamma (a+k) \Gamma (b+2 k) \Gamma (d+k) \Gamma (2 a-b-2 d+1)}{\Gamma (b) \Gamma (d) \Gamma (k+1) \Gamma \left(a-b-d+\frac{1}{2}\right) \Gamma (2 a-b+2 k+1) \Gamma (a-d+k+1)}.$$ Let $F(n,k) = \textsf{Seed3}(a+A n,b+Bn,d+Dn, k+k_0+Kn)$, exponent of $F(n,0)$ easily seen to be\footnote{here we agree $0^0 = 1$} \begin{multline*}B^{-B} 2^{2 D} D^{-D} K^{-K} (A-B)^{A-B} (A+K)^{A+K} (B+2 K)^{B+2 K} \\ \times (D+K)^{D+K} (2 A-B-2 D)^{2 A-B-2 D} (A-B-D)^{-A+B+D} (2 A-B+2 K)^{-2 A+B-2 K} (A-D+K)^{-A+D-K}.\end{multline*}
It depends on accessory parameters $A,B,D,K$, and is independent of $a,b,d,k$. If we let $A=1,B=-1,D=1,K=-1$, then we obtain $16/27$, our $F(n,0)$ is $$-\frac{4^{d+n} \left(-a-2 k_0+n\right) \Gamma \left(a+k_0\right) \Gamma \left(d+k_0\right) \Gamma \left(a-b+2 n+\frac{1}{2}\right) \Gamma \left(b-3 n+2 k_0\right) \Gamma (2 a-b-2 d+n+1)}{\Gamma (b-n) \Gamma (d+n) \Gamma \left(-n+k_0+1\right) \Gamma \left(a-b-d+n+\frac{1}{2}\right) \Gamma \left(2 a-b+n+2 k_0+1\right) \Gamma \left(a-d-n+k_0+1\right)}$$
now pochhammer part $\frac{(a_1)_n\cdots (a_m)_n}{(b_1)_n \cdots (b_m)_n}$ (recall we required $a_i,b_i \in [0,1)$) will depend on $(a,b,c,k_0) \in \mathbb{Q}/\mathbb{Z}$. For example, for $(a,b,c,k_0) \equiv (1/3,1/3,1,1/3)\pmod{\mathbb{Z}}$, the pochhammer part will be $\frac{(1/4)_n (1/3)_n (2/3)_n (3/4)_n}{(1/2)_n (1)_n^3}$. Therefore we conclude, if we want to prove a conjectural identity of the form $$\sum_{n\geq 0} (\frac{16}{27})^n \frac{(1/4)_n (1/3)_n (2/3)_n (3/4)_n}{(1/2)_n (1)_n^3} R(n) \in \frac{\text{CMZVs of some level}}{\pi^r},$$ we should apply Proposition \ref{WZ_prop} on $F(n,k) = \textsf{Seed3}(a + n,c -n,d+n, k+k_0-n)$ with $(a,b,c,k_0) \equiv (1/3,1/3,1,1/3)\pmod{\mathbb{Z}}$. 

To summarize our execution:
\begin{tcolorbox}[colback=gray!5!white,colframe=gray!75!black,title= Searching suitable WZ-candidate for a hypergeometric identity]
	Input: a collection of WZ-seeds, $z$ and $a_i,b_i$ in Equation \ref{to_be_proved_identity}. \par
	Output: A WZ-seed \textsf{Seed}, a choice of $A,B,\cdots, K \in \mathbb{Z}$, a choice of $a,b,\cdots,k_0 \in \mathbb{Q}/\mathbb{Z}$ such that for $F(n,k) = \textsf{Seed}(a+An,b+Bn,k+k_0+Kn)$, $G(n,0)$ is up to rational function, $z^n \dfrac{(a_1)_n\cdots (a_m)_n}{(b_1)_n \cdots (b_m)_n}$. 
	\tcblower
\begin{enumerate}
	\item Choose a small range of integers, say from $-N$ to $N$ with $N\geq 1$. Choose a small positive integer, say $l\geq 1$.
	\item For each $\textsf{Seed}(a,b,\cdots,k)$ in the collection, do
	\begin{enumerate}
		\item For each possible combination of $A,B,\cdots,K \in [-N,N]$, compute the exponent of $\textsf{Seed}(a+An,b+Bn,\cdots,k+k_0+Kn)$, check if any one of them equals $z$.
		\item If no, skip to (3). If yes, for each possible combination of $a,b,\cdots,k_0 \in \{1/l,2/l,\cdots,1\}$, compute the pochhammer part of $\textsf{Seed}(a+An,b+Bn,\cdots,k+k_0+Kn)$, check if any one of them equals $\frac{(a_1)_n\cdots (a_m)_n}{(b_1)_n \cdots (b_m)_n}$
		\item If no, skip to (3). If yes, return $A,B\cdots,K$ and $(a,b,\cdots,k_0)\in \mathbb{Q}/\mathbb{Z}$
	\end{enumerate}
	\item Go to next seed and do (2). 
\end{enumerate}
\end{tcolorbox}
A Mathematica implementation for above procedure can be found in \url{https://www.researchgate.net/publication/376812577}. All WZ-pairs used in this article are generated from it. 
For our previous two examples, $F(n,k)$ in Example \ref{ex_pi-4_1} comes from $\textsf{Dougall7F6}$ with $A=B=C=D=E=-1, K=2$ and $a=b=c=d=e=1/2,k_0=1$. $F(n,k)$ in Example \ref{ex_pi-4_2} comes from $\textsf{Dougall7F6}$ with $A=1,B=C=D=E=0, K=1$ and $a=b=c=d=e=1/2,k_0=1$. 

~\\[0.01in]
\textit{With no exception known, once the above searching finds a match, conjectural identity \ref{to_be_proved_identity} is destined to be proved}. There are however, sometimes plenty of loose ends to be tied. First and foremost, we only determined $a,b,\cdots,k_0$ modulo $\mathbb{Z}$, more trial and error is required to pin-down an exact choice\footnote{in some rare cases, this seems not possible, then $G(n,0)$ and to-be-proved identity \ref{to_be_proved_identity} will differ by a telescoping series, see Examples \ref{384pi-2_ex} and \ref{12pi-2_ex}.} so that $G(n,0)$ gives exactly the desired identity (\ref{to_be_proved_identity}). Secondly, in the statement of Proposition \ref{WZ_prop}, any of the three terms $$\sum_{k\geq 0} F(0,k), \qquad \sum_{n\geq 0} g(n), \qquad \lim_{n\to\infty} \sum_{k\geq 0} F(n,k),$$
might be a source of nuisance:
\begin{itemize}
\item $g(n)$ might not be zero, so we still have to evaluate $\sum_{n\geq 0}g(n)$. Fortunately this $g(n)$ will always be of lesser complexity than $G(n,0)$, so we can repeat the above search to $g(n)$ again\footnote{for all cases in this article, $g(n)$ becomes so much simpler that they're already known in the literature, as in previous two examples. This spares us from performing this recursive step.} looking for another $(F',G')$ such that $G'(n,0) = g(n)$. 
\item It might not be immediately clear what is $\sum_{k\geq 0} F(0,k)$. Our way to overcome this is hypergeometric transformation formulas. (e.g. Examples \ref{ex_CMZVtrans}, \ref{pi-2_parameters})
\item The most intractable term is $\lim_{n\to\infty} \sum_{k\geq 0} F(n,k)$, we will investigate when it will be zero below. An entire section is devoted to examples in which it is non-zero. 
\end{itemize}
The variables $a,b,\cdots,k_0$ can be kept as indeterminate, they are free variables in the final formula of $\sum_{n\geq 0}G(n,0)$. For a WZ-seed, the number of free variables in final formula is the same as the number of arguments in the seed. In previous two examples, we collapsed five (out of six) free variables, otherwise the formula would be way to long to display. Example \ref{gauss_ex1} below retains all free variables. In the paper, we translate all free variables $a,b,\cdots$ to $0$. \par Advantage of retaining these free variables is already evident in previous two examples: setting $a=0$ straight away would not give a viable equality. Moreover, they can be used to prove identity (\ref{to_be_proved_identity}) twisted with harmonic numbers. 
~\\[0.01in]

We mention an limitation of our method: let $\mathcal{G} \subset \mathbb{Q}^\times$ be the subgroup generated by $\{n^n | n\in \mathbb{Z}\}$. From the lists of seeds, only \textsf{Seed1}, \textsf{Seed4}, \textsf{Seed10} can generate series whose exponent of $G(n,0)$ is not in $\mathcal{G}$. Therefore many (known or conjectural) series, for instances, 
$$\begin{aligned}\sum_{n\geq 0} (\frac{-1}{80^3})^n \frac{(1/2)_n (1/6)_n (5/6)_n}{(1)_n^3} (5418n+263) &= \frac{640\sqrt{15}}{3\pi},\\
\sum_{n\geq 0} (\frac{1}{7^{4}})^n \frac{(1/2)_n (1/8)_n (3/8)_n (5/8)_n (7/8)_n}{(1)_n^5} (1920n^2+304n+15) &\stackrel{?}{=} \frac{56\sqrt{7}}{\pi^2}\end{aligned}$$
do not succumb under this technique.

\subsection{Vanishing of limit $\lim_{n\to\infty} \sum_{k\geq 0} F(n,k)$}
In this section, we give a convenient sufficient criterion on when $\lim_{n\to\infty} \sum_{k\geq 0} F(n,k)$ vanishes. Let $$F(n,k) = A^n B^k \frac{\prod_{i=1}^m\Gamma(a_i + b_i n + c_i k)}{\prod_{i=1}^r\Gamma(d_i + e_i n + f_i k)},$$
here $a_i,b_i,c_i,d_i,e_i,f_i \in \mathbb{R}, A, B\in \mathbb{C}$, moreover, we also stipulate $\sum_i b_i = \sum_i e_i$ and $\sum_i c_i = \sum_i f_i$. All examples of $F(n,k)$ in this article are in this form. We do not require it has a hypergeometric WZ-mate $G(n,k)$. 

We lighten our notation with the shorthand $p(x) = \begin{cases}|x|^x \quad x\neq 0 \\ 1 \quad x=0\end{cases}$. To such $F$, we define $\mathcal{E}_F(x)$, which is a positive-valued function on $[0,\infty)$, 
$$\mathcal{E}_F(x) = |A|\times |B|^x\times \frac{\prod_{i=1}^m p(b_i + c_i x)}{\prod_{i=1}^r p(e_1 + f_i x)}.$$
It is easy to check, using Stirling's formula, for each $x\geq 0$,
$$|F(n,xn)| = (\mathcal{E}_F(x))^n \times O(n^C), \qquad C\in \mathbb{R}, \qquad n\to \infty.$$
Here $C$ also depends on $a_i,d_i$ in $F(n,k)$. It also follows, from this observation, that $\mathcal{E}_F(x)$ is independent of the representation in terms of gamma function we used to write $F$. Moreover, $C$ and the constant in $O$-term are uniformly bounded for $x\geq 0$. 

\begin{proposition}\label{vanishing_prop}
$\lim_{n\to\infty} \sum_{k\geq 0} F(n,k)$ is $0$ if $$\sup\{\mathcal{E}_F(x) | x\geq 0 \} < 1.$$
\end{proposition}
\begin{proof}
Let $S$ be the suprenum, we pick $\varepsilon > 0$ such that $S+\varepsilon < 1$. For given $r\in \mathbb{N}$, we have $$\left|\sum_{rn \leq k < (r+1)n} F(n,k) \right| < n \times (S+\varepsilon)^n $$ for sufficiently large $n$, so tends to $0$ as $n\to \infty$. Hence $$\sum_{k\geq 0} F(n,k) = \sum_{r\geq 0} \sum_{rn \leq k < (r+1)n} F(n,k)$$
also tends to $0$ by dominated convergence.
\end{proof}

We can now quickly check, in Examples \ref{ex_pi-4_1} and \ref{ex_pi-4_2}, the term $\lim_{n\to\infty} \sum_{k\geq 0} F(n,k) = 0$. For $F$ in former example, we have
$$\mathcal{E}_F(x) = \frac{| x-1| ^{1-x} | x+1| ^{5 x+5} | x+2| ^{-5 x-10} | x+4| ^{x+4}}{1024},$$
it satisfies $\sup\{\mathcal{E}_F(x) | x\geq 0 \} < 1$ is an easy calculus exercise. 
\begin{figure}[h]
	\centering
	\includegraphics[width=0.6\textwidth]{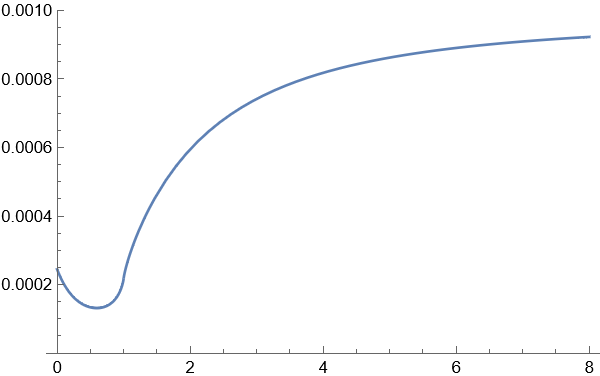}
	\caption{Plot of $\mathcal{E}_F(x)$, with $F$ in Example \ref{ex_pi-4_1}}
	\label{plot}
\end{figure}
Similarly one can check $\mathcal{E}_F(x)$ in Example \ref{ex_pi-4_2} also meets this sufficient condition to guarantee vanishing of the limit. All examples in next section satisfy Proposition \ref{vanishing_prop}, so $\lim_{n\to\infty} \sum_{k\geq 0} F(n,k) = 0$ for them.

\subsection{Euler sums and CMZVs}
We assume familiarity with colored multiple zeta values (CMZVs). For more background on this subsection, we refer readers to author's previous article \cite{au2022multiple}. 

The series
$$\Li_{s_1,\cdots,s_k}(a_1,\cdots,a_k) = \sum_{n_1>\cdots>n_k\geq 1}\frac{a_1^{n_1}\cdots a_k^{n_k}}{n_1^{s_1} \cdots n_k^{s_k}}$$
is known as \textit{multiple polylogarithm}. Here $k$ is called the \textit{depth} and $s_1+\cdots+s_k$ is called the \textit{weight}. 

When $a_i$ are $N$-th roots of unity, $s_i$ are positive integers and $(a_1, s_1) \neq (1,1)$, $\Li_{s_1,\cdots,s_k}(a_1,\cdots,a_k)$ is called a \textit{colored multiple zeta values} (CMZV) of weight $s_1+\cdots+s_k$ and \textit{level} $N$. We say a complex number is CMZV of weight $w$ and level $N$ if it's in $\mathbb{Q}(e^{2\pi i /N})$-span of such numbers\footnote{in other literature, sometimes one takes $\mathbb{Q}$-span instead of $\mathbb{Q}(e^{2\pi i /N})$-span}, which we denote by $\CMZV{N}{w}$. The special case when all $a_i = 1$ is the well-known \textit{multiple zeta value}.

	Fix a positive integer $N$, $0<\gamma_i \leq 1$ , we denote
$$H(\vec{s},\vec{\gamma},N) = \sum_{\substack{n_1+\gamma_1 > n_2+\gamma_2 > \cdots > n_k+\gamma_k > 0 \\ n_1<N}} \frac{1}{(n_ 1+\gamma_1)^{s_ 1}\cdots (n_k+\gamma_k)^{s_k}},$$
here $n_i$ in sum need to be integers; with $\vec{s} = (s_1,\cdots,s_k)\in 
\mathbb{N}^k,\vec{\gamma} = (\gamma_1,\cdots,\gamma_k) \in \mathbb{Q}^k$. 



For $0<\gamma\leq 1$, we have easy-to-prove series expansion (around $a=0$):
\begin{equation}\label{pochhammerexpansion}(\gamma+a)_n = (\gamma)_n \left(1+ \sum_{k\geq 1} H(\{1\}_k, \{\gamma\}_k, n) a^k \right).\end{equation}

The following proposition was proved in \cite{au2022multiple}. \par

	\begin{proposition}\label{CMZVsum1}
	Let $0<\gamma_0\leq 1$, $K_i$ positive integers, $\vec{\gamma_i}$ and $\vec{a_i}$ rational vectors such that each component is $0< \cdot \leq 1$. Let $N$ be an integer such that all $N\vec{\gamma_i}, N\vec{a_i}, N\gamma_0, Na_0, N/K_i$ are in $\mathbb{Z}$. The following series, regularized (in the sense defined in \cite{au2022multiple}), is a level $N$ (possibly inhomogeneous) CMZV
	$$\sum_{n_0>0} \frac{e^{2\pi i  a_0 n_0}}{(n_0+\gamma_0)^{s_0}} H(\vec{s_1},\vec{\gamma_1},\vec{a_1},K_1 n_0) \cdots H(\vec{s_m},\vec{\gamma_m},\vec{a_m},K_m n_0) \in \sum_{1\leq i\leq w} \CMZV{N}{i},$$
	where $w = |\vec{s_1}| + \cdots + |\vec{s_i}| + s$. If $\gamma_0 = 1/N$, then it is a homogeneous CMZV of weight $w$. 
\end{proposition}

Converting sums into CMZVs is desirable, largely because we now understand (at least conjecturally, for small weight and level) all algebraic relations between them, so any purported algebraic relation between them can be trivially verified. \\[0.01in]

As a corollary of Proposition \ref{CMZVsum1} and equation (\ref{pochhammerexpansion}), for the hypergeometric series
$$\pFq{r+1}{r}{A_1 + a_1,A_2+a_2,\cdots A_{r+1}+a_{r+1}}{B_1 + b_1,B_2+b_2,\cdots B_{r}+b_{r}}{z},\qquad A_i,B_i \in \mathbb{Q},$$ 
if $z$ is a root of unity and \begin{equation}\label{int_diff_cond}\text{difference between a permutation of } (A_1,\cdots,A_{r+1}) \text{ and } (1,B_1,\cdots,B_{r}) \in \mathbb{Z}^{r+1}.\end{equation} Then, when the function are expanded in power series in variables $a_1.\cdots,a_{r+1},b_1,\cdots,b_r$ around origin, each coefficient will be (possibly inhomogeneous in weight) CMZVs of some level. 
\\[0.01in]

For later usage, we assemble here some notations and formulas of very-well-poised hypergeometric series. Denote
$$V(a_1;a_2,\cdots,a_r) = \pFq{r+1}{r}{a_1,1+\frac{a_1}{2},a_2,a_3,\cdots,a_r}{\frac{a_1}{2},1+a_1-a_2,1+a_1-a_3,\cdots,1+a_1-a_r}{(-1)^r},$$
then for $_7F_6$, we have transformations (\cite{bailey1935generalized}):
\begin{align}\label{VWP7F6trans} \begin{split} & V(a;c,d,e,f,g)\\  \quad &= \frac{\Gamma (a-c+1) \Gamma (a-d+1) \Gamma (2 a-e-f-g+2) \Gamma (2 a-c-d-e-f-g+2)}{\Gamma (a+1) \Gamma (a-c-d+1) \Gamma (2 a-c-e-f-g+2) \Gamma (2 a-d-e-f-g+2)}\\ & \times V(1 + 2 a - e - f - g; c, d, 1 + a - f - g, 1 + a - e - g, 1 + a - e - f) \\
		&= \frac{\Gamma (a-c+1) \Gamma (a-d+1) \Gamma (a-e+1) \Gamma (a-f+1) \Gamma (3 a-c-d-e-f-2 g+3) \Gamma (2 a-c-d-e-f-g+2)}{\Gamma (a+1) \Gamma (g) \Gamma (2 a-c-d-e-g+2) \Gamma (2 a-c-d-f-g+2) \Gamma (2 a-c-e-f-g+2) \Gamma (2 a-d-e-f-g+2)} \\ & \times V(3 a-c-d-e-f-2 g+2;a-c-g+1,a-d-g+1,a-e-g+1,a-f-g+1,2 a-c-d-e-f-g+2).\end{split}\end{align}

This and many other transformations enable us to convert a $_{r+1}F_r$ that does not satisfy condition (\ref{int_diff_cond}) to one that possibly does. 

	\section{Basic examples}
	
	Some examples below will involve either heavy computation or formulas too long to be displayed. Readers interested in computational details and full expressions could refer to a Mathematica file at  \url{https://www.researchgate.net/publication/376812577}. 
	
	\begin{romexample}\label{gauss_ex1}
Let $p_n = 42n^2-23n+3$, $a_n = 2^{-6n} \dfrac{(1/2)_n (1)_n^3}{(1/4)_n^2 (3/4)_n^2}$, then Z.W. Sun conjectured\footnote{\url{https://mathoverflow.net/questions/456443}}
	$$\begin{aligned}\sum_{n\geq 1} a_n \frac{p_n}{n^3(2n-1)} &= \frac{\pi^2}{2}, \\
	\sum_{n\geq 1} a_n \frac{p_n (H_{2n-1}-H_{n-1}) - \frac{196n^2-100n+13}{6(2n-1)}}{n^3(2n-1)} &= \frac{\pi^2 \log 2}{3} - \frac{7\zeta(3)}{6},\\
	\sum_{n\geq 1} a_n \frac{p_n (H_{4n-1}-H_{2n-1}) - \frac{28n^2-76n+19}{12(2n-1)}}{n^3(2n-1)} &= \frac{\pi^2 \log 2}{6} + \frac{35\zeta(3)}{12},\\
	\sum_{n\geq 1} a_n \frac{p_n (4H^{(2)}_{4n-1}-H^{(2)}_{2n-1}-2H^{(2)}_{n-1}) - 6-1/n}{n^3(2n-1)} &= \frac{5\pi^4}{24}, \\
	\sum_{n\geq 1} a_n \frac{p_n (4H^{(2)}_{4n-1}-5H^{(2)}_{2n-1}-3H^{(2)}_{n-1}) + \frac{32n(3n-1)}{(2n-1)^2}}{n^3(2n-1)} &= \frac{\pi^4}{6}.
		\end{aligned}$$
	
	Searching for exponent $1/64$ and hypergeometric part $\frac{(1/2)_n (1)_n^3}{(1/4)_n^2 (3/4)_n^2}$ through our lists of WZ-seeds returns a match under the seed \textsf{Gauss2F1}, and we are inspired to take $$F(n,k) = \frac{\Gamma (-a+c+n+1) \Gamma \left(a+k+n+\frac{1}{2}\right) \Gamma (-b+c+n+1) \Gamma \left(b+k+n+\frac{1}{2}\right)}{\Gamma (a-n) \Gamma (b-n) \Gamma \left(k+2 n+\frac{3}{2}\right) \Gamma \left(c+k+2 n+\frac{3}{2}\right) \Gamma (-a-b+c+2 n+1)}.$$ One finds the corresponding $G(n,k)$, $g(n) = 0$. $\lim_{n\to \infty} \sum_{k\geq 0}F(n,k) = 0$ because $\mathcal{E}_F(x)$ meets condition of Proposition \ref{vanishing_prop}. Now $$\sum_{k\geq 0} F(0,k+d) = \sum_{n\geq 1} G(n-1,d)$$ becomes
\begin{multline*}\sum_{k\geq 0}\frac{32 \left(a+d+\frac{1}{2}\right)_k \left(b+d+\frac{1}{2}\right)_k}{(2 d+2 k+1) (2 c+2 d+2 k+1) \left(d+\frac{1}{2}\right)_k \left(c+d+\frac{1}{2}\right)_k} \\ = \sum_{n\geq 1}\frac{(1-a)_{n-1} (1-b)_{n-1} (-a+c+1)_{n-1} \left(a+d+\frac{1}{2}\right)_{n-1} (-b+c+1)_{n-1} \left(b+d+\frac{1}{2}\right)_{n-1}  \times P}{\left(d+\frac{1}{2}\right)_{2 n} \left(c+d+\frac{1}{2}\right)_{2 n} (-a-b+c+1)_{2 n}},\end{multline*}
	with $P \in \mathbb{Z}[a,b,c,d,n]$ is a long polynomial that pops up when computing WZ-mate, its full-form can be found in the Mathematica file linked. When $a=b=c=d=0$, the formula reduces to the first conjecture above, those remaining four involving harmonic numbers can be obtained by comparing coefficients at $(a,b,c,d) = (0,0,0,0)$, with coefficients of LHS being level 2 CMZVs by Proposition \ref{CMZVsum1}. 
	
	\end{romexample}
	
	\begin{romexample}\label{ex_pi4_1}J. Guillera conjectures (\cite{guillera2003new}) 
		$$\sum_{n\geq 1} (\frac{4}{27})^n \frac{(1)_n^7}{(1/2)_n^5 (1/3)_n (2/3)_n} \frac{92n^3-84n^2+21n-3}{n^7} = 8\pi^4.$$ 
		Searching for exponent $4/27$ and  $\frac{(1)_n^7}{(1/2)_n^5 (1/3)_n (2/3)_n}$ yields a match under \textsf{Seed7}. We are inspired to take
		$$F(n,k) = \frac{n \left(2 k+2 n+\frac{3}{2}\right) \Gamma (n+1)^3 \Gamma (k+n+1)^3 \Gamma (2 k+n+1) \Gamma (k+2 n+1)}{\Gamma \left(k+\frac{3}{2}\right) \Gamma \left(n+\frac{1}{2}\right)^3 \Gamma \left(k+n+\frac{3}{2}\right)^3 \Gamma (2 k+3 n+3)}.$$
		
		After finding $G(n,k), g(n)$, 
		$$\sum_{k\geq 0} F(a,k) + \sum_{n\geq 0} g(n+a) = \sum_{n\geq 0} G(n+a,0)$$ becomes
		
		\begin{multline*}\sum_{k\geq 0} \frac{8 (4 a+4 k+3) (1+a)_k^3 (1+a)_{2 k} (2 a+1)_k}{(2 k+1) (2 a+2 k+1)^3 (3 a+2 k+1) (3 a+2 k+2) (1/2)_k (1/2+a)_k^3 (1+3a)_{2k}}  \\ - \sum_{n\geq 0} \frac{\sqrt{\pi } 2^{-2 a-2 n+1} (3 a+3 n+2) \Gamma(a+\frac{1}{2})^3 \Gamma (3 a+1) (1+a)_n^3}{a (2 a+2 n+1)^3 \Gamma (a+1)^4 \Gamma (2 a+1) (1/2+a)_n^3} \\  = - \sum_{n\geq 1} \frac{(92 a^3+276 a^2 n-84 a^2+276 a n^2-168 a n+27 a+92 n^3-84 n^2+27 n-3) (1+a)_n^7 (2 a+1)_{2 n}}{16 a (a+n)^7 (1/2+a)_n^6 (3 a+1)_{3 n}}.\end{multline*}
		
		Comparing coefficient of $a^{-1}$ gives (again the term $\sum_{k\geq 0} \cdots$ does not contribute):
		$$\sum_{n\geq 0} \frac{\pi ^2 2^{1-2 n} (3 n+2)(1)_n^3}{(2 n+1)^3 (1/2)_n^3} = \sum_{n\geq 1} \frac{\left(92 n^3-84 n^2+27 n-3\right) (1)_n^7 (1)_{2 n}}{16 n^7 (1/2)_n^6 (1)_{3 n}}.$$
		Invoking the known WZ-type result $$\sum_{n\geq 0} \frac{2^{-2 n} (3 n+2)(1)_n^3}{(2 n+1)^3 (1/2)_n^3} = \frac{\pi^2}{4} $$ completes the proof.
	\end{romexample}

	\begin{romexample}\label{ex_CMZVtrans}In this example, we look at an example for which $\sum_{k\geq 0} F(0,k)$ is not in an immediately useable form, some further manipulation is thus needed. We motivate ourselves with the two conjectures of Z.W. Sun \footnote{\url{https://mathoverflow.net/questions/456444}}, here $L_{-3}(s) = \sum_{n\geq 0} \left((3n+1)^{-s} - (3n+2)^{-s}\right)$, 
		$$\sum_{n\geq 1} (-\frac{1}{48})^n \frac{(\frac{1}{2})_n (1)_n^3}{(\frac{1}{4})_n (\frac{1}{3})_n (\frac{2}{3})_n (\frac{3}{4})_n} \frac{42n^2-27n+4}{n^3 (2n-1)} = -6 L_{-3}(2),$$
		$$\sum_{n\geq 1} (-\frac{1}{48})^n \frac{(\frac{1}{2})_n (1)_n^3}{(\frac{1}{4})_n (\frac{1}{3})_n (\frac{2}{3})_n (\frac{3}{4})_n} \frac{(42n^2-27n+4)(H_{3n-1}-H_{n-1})-10n+3}{n^3 (2n-1)} =  -\frac{8\pi^3}{27\sqrt{3}}.$$
		
		Apply Proposition \ref{WZ_prop} to 
		$$F(n,k) = \frac{3^{-n} \Gamma (3 a-2 b+2 n+2) (3 a+6 c+6 k+6 n+4) \Gamma (a+c+k+n+1) \Gamma (b+3 c+3 k+2 n+1)}{\Gamma (b-n) \Gamma \left(a-b+n+\frac{2}{3}\right) \Gamma \left(c+k+n+\frac{4}{3}\right) \Gamma (3 a-b+3 c+3 k+4 n+4)},$$
		one quickly checks $\lim_{n\to\infty} \sum_{k\geq 0} F(n,k) = 0$ and $g(n) = 0$. So $\sum_{n\geq 1}G(n-1,0) = \sum_{k\geq 0} F(0,k)$, that is, 
		\begin{multline*}(3 a-2 b+1)\sum_{k\geq 0}\frac{(3 a+6 c+6 k+4) (a+c+1)_k (b+3 c+1)_{3 k}}{(3 c+3 k+1) (c+\frac{1}{3})_k (3 a-b+3 c+1)_{3 k+3}} \\ = \sum_{n\geq 1} \frac{(-1)^{n+1} 3^{-n-1} (1-b)_{n-1} (3 a-2 b+1)_{2 n} (a+c+1)_{n-1} (b+3 c+1)_{2 n-2}  \times P}{(3 a-2 b+2 n) (c+\frac{1}{3})_n (a-b+\frac{2}{3})_n (3 a-b+3 c+1)_{4 n}},\end{multline*}
		with $P\in \mathbb{Z}[a,b,c,n]$ whose full-form can be in the file mentioned at start of this section. \par
		 When $a=b=c=0$, RHS becomes the first conjecture of Sun above, and the other one follows by taking linear combination of coefficients. However, it's not immediately clear that LHS's coefficients are CMZVs, as Proposition \ref{CMZVsum1} and \ref{int_diff_cond} no longer applies. We claim LHS's coefficient are indeed level 3 CMZVs. \par
		 
		 To see this, write LHS as $$\frac{(1+3a-2b)(4+3a+6c)}{(1 + 3 c) (1 + 3 a - b + 3 c) (2 + 3 a - b + 3 c) (3 + 3 a - b + 3 c)} \pFq{5}{4}{1,a+c+1,\frac{b}{3}+c+\frac{1}{3},\frac{b}{3}+c+\frac{2}{3},\frac{b}{3}+c+1}{c+\frac{4}{3},a-\frac{b}{3}+c+\frac{4}{3},a-\frac{b}{3}+c+\frac{5}{3},a-\frac{b}{3}+c+2}{1}.$$
		 This $_5F_4$ is actually a very-well-poised $_7F_6$, $$V(a+2 c+\frac{4}{3};1,a+c+1,\frac{b}{3}+c+\frac{1}{3},\frac{b}{3}+c+\frac{2}{3},\frac{b}{3}+c+1).$$
		 By transformation \ref{VWP7F6trans}, this is
		 \begin{multline*}\frac{\Gamma \left(a-b+\frac{2}{3}\right) \Gamma \left(a+c-\frac{b}{3}+\frac{4}{3}\right) \Gamma \left(a+c-\frac{b}{3}+2\right) \Gamma \left(a+2 c-\frac{b}{3}+2\right)}{\Gamma \left(a-\frac{2 b}{3}+1\right) \Gamma \left(a+2 c+\frac{7}{3}\right) \Gamma \left(a+c-\frac{2 b}{3}+1\right) \Gamma \left(a+c-\frac{2 b}{3}+\frac{5}{3}\right)} \\ \times V(a-\frac{b}{3}+2 c+1;\frac{b}{3}+c+\frac{1}{3},\frac{b}{3}+c+1,\frac{2}{3}-\frac{b}{3},a-\frac{b}{3}+c+\frac{2}{3},c+\frac{1}{3}).\end{multline*}
		 
		 Therefore our original sum $\sum_{k\geq 0}\cdots$ equals
		 \begin{multline*}(1+3a-2b)\frac{\Gamma (a-b+\frac{2}{3}) \Gamma (a+c-\frac{b}{3}+\frac{1}{3}) \Gamma (a+c-\frac{b}{3}+1) \Gamma \left(a+2 c-\frac{b}{3}+1\right)}{\Gamma (a-\frac{2 b}{3}+1) \Gamma (a+2 c+\frac{1}{3}) \Gamma (a+c-\frac{2 b}{3}+\frac{2}{3}) \Gamma (a+c-\frac{2 b}{3}+1)} \\ \times \sum_{k\geq 0} \frac{(3+3a-b+6c+6k)}{(3 c+3 k+1) (3 a+6 c+3 k+1) (3 a-2 b+3 c+3 k+2) (3 a-b+3 c+3 k+2)} \\ \times \frac{(2/3-b/3)_k (1/3+b/3+c)_k (1+b/3+c)_k (1+a-b/3+2c)_k}{(1)_k (2/3+a-2b/3+c)_k (1+a-2b/3+c)_k (1/3+a+2c)_k}. \end{multline*}
		 
		 in this form, equation (\ref{int_diff_cond}) now holds and we conclude the coefficients of $a^ib^jc^k$ are indeed level 3 CMZVs (of weight $i+j+k+2$). 
	\end{romexample}

\begin{romexample}\label{pi-2_parameters}
Let $a_n = (-\frac{1}{2^{10}})^n \frac{(1/2)_n^5}{(1)_n^5}$, $p_n = 820n^2+180n+13$. Consider the following twists by harmonic numbers of the famous $\sum_{n\geq 0} a_n p_n = \frac{128}{\pi^2}$:
$$\begin{aligned}
	\sum_{n\geq 0} a_n (p_n (H_{2n}-H_n) + 164n+18) &= \frac{256\log 2}{\pi^2},\\
	\sum_{n\geq 0} a_n (p_n (11 H_{2 n}^{(2)}-3 H_n^{(2)}) + 43) &= \frac{128}{3},\\
	\sum_{n\geq 0} a_n (p_n (9 H_{2 n}^{(3)}-H_n^{(3)}) + \frac{125}{2n+1}) &= \frac{2^{10} \zeta(3)}{\pi^2},\\
	\sum_{n\geq 0} a_n (p_n (49 H_{2 n}^{(4)}-3 H_n^{(4)}) - \frac{195}{(2n+1)^2}) &= -\frac{896}{45}\pi^2, \\
	\sum_{n\geq 0} a_n (p_n (351 H_{2 n}^{(5)}-11 H_n^{(5)})+\frac{275}{(2 n+1)^3}) &= 2^9\left(\frac{85\zeta(5)}{\pi^2}-7\zeta(3)\right).
\end{aligned}$$
All were conjectured by Z.W. Sun (\cite{sun2023new}, \cite{sun2022conjectures}). First two are proved in \cite{wei2023some2}, third in \cite{wei2023some3}. fourth in \cite{hou2023taylor}. Their methods, which are all based on a WZ-type formula produced by Chu and Zhang (\cite{chu2014accelerating}), given enough computations, can actually prove all of above. We shall re-establish them by extracting coefficients from a $5$-parameter deformation of $\sum_{n\geq 0} a_n p_n = \frac{128}{\pi^2}$:
\begin{multline*}
A\times \biggl[1+(2 a-d+e) (2 a-b-c-d) \sum_{k\geq 0} \frac{\splitfrac{ \left(a-d+\frac{1}{2}\right)_{k+1} \left(a-b-d+\frac{1}{2}\right)_{k+1}}{ \left(a-c-d+\frac{1}{2}\right)_{k+1} (4 a-b-c-2 d+e+2 k+2) (4 a-b-c-2 d+e+1)_k}}{\splitfrac{(k+1) (1)_k (2 a-d+e+k+1) (2 a-b-c-d+k+1) \left(3 a-b-d+e+\frac{1}{2}\right)_{k+1}}{ \left(3 a-c-d+e+\frac{1}{2}\right)_{k+1} \left(3 a-b-c-d+e+\frac{1}{2}\right)_{k+1}}} \biggl] \\ = \sum_{n\geq 0} \frac{\splitfrac{(-1)^n P(n+a,e) \left(a-b+\frac{1}{2}\right)_n \left(a-c+\frac{1}{2}\right)_n \left(a-d+\frac{1}{2}\right)_n \left(a+e+\frac{1}{2}\right)_n\times }{ \left(a-b-c+\frac{1}{2}\right)_n \left(a-b-d+\frac{1}{2}\right)_n \left(a+b+e+\frac{1}{2}\right)_n \left(a-c-d+\frac{1}{2}\right)_n \left(a+c+e+\frac{1}{2}\right)_n \left(a+d+e+\frac{1}{2}\right)_n}}{(2 a+e+1)_{2 n+1} (2 a-b+e+1)_{2 n+1} (2 a-c+e+1)_{2 n+1} (2 a-d+e+1)_{2 n+1} (2 a-b-c-d+1)_{2 n+1}}.
\end{multline*}
where $P(n,k) \in \mathbb{Z}[b,c,d,n,k]$ is a long polynomial that has $608$ monomial terms\footnote{readers interested in full expression should consult the Mathematica notebook attached} and $$A = 128\frac{\Gamma (2 a+e+1) \Gamma (2 a-b+e+1) \Gamma (2 a-c+e+1) \Gamma (4 a-b-c-2 d+e+1)}{\Gamma \left(a+d+e+\frac{1}{2}\right) \Gamma \left(3 a-b-d+e+\frac{1}{2}\right) \Gamma \left(3 a-c-d+e+\frac{1}{2}\right) \Gamma \left(3 a-b-c-d+e+\frac{1}{2}\right)}.$$
Above reduces to $\sum_{n\geq 0}a_np_n = \frac{128}{\pi^2}$ when all five parameters $a,b,c,d,e$ are $0$, the $\pi^{-2}$ is explained by the gamma product term $A$. In general, comparing coefficient of degree $w$ on LHS will yield $$\pi^{-2}\times (\text{CMZVs of level 2 and weight }w)$$ by equation \ref{int_diff_cond}. The above formula follows from Proposition \ref{WZ_prop} with \footnote{after searching for exponent $-1/2^{10}$ and pochhammer part $(1/2)_n^5/(1)_n^5$} $$F(n,k) = \frac{\splitfrac{(4 k+6 n+1) \Gamma \left(k+n+\frac{1}{2}\right) \Gamma \left(-b-c+n+\frac{1}{2}\right) \Gamma \left(-b-d+n+\frac{1}{2}\right) \Gamma \left(b+k+n+\frac{1}{2}\right)}{ \Gamma \left(-c-d+n+\frac{1}{2}\right) \Gamma \left(c+k+n+\frac{1}{2}\right) \Gamma \left(d+k+n+\frac{1}{2}\right)}}{\splitfrac{\Gamma \left(b-n+\frac{1}{2}\right) \Gamma \left(c-n+\frac{1}{2}\right) \Gamma \left(d-n+\frac{1}{2}\right)}{ \Gamma (k+2 n+1) \Gamma (-b+k+2 n+1) \Gamma (-c+k+2 n+1) \Gamma (-d+k+2 n+1) \Gamma (-b-c-d+2 n)}},$$
the identity is $$\sum_{k\geq 0} F(a,k+e) = \sum_{n\geq 0} G(n+a,e)$$
together with some transformation of very-well-poised $_7F_6$ on $\sum_{k\geq 0} F(a,k+e)$ such that it's coefficient at $(a,b,c,d,e) = (0,0,0,0,0)$ obviously becomes CMZVs, a process similar to that in Example \ref{ex_CMZVtrans}. 
\end{romexample}

\begin{romexample}\label{384pi-2_ex}Consider the following conjectural identity that appears in \cite{guillera2016bilateral}, 
	$$\sum_{n\geq 0} a_n (1930n^2+549n+45) = \frac{384}{\pi^2}, \qquad a_n = (-\frac{3^6}{2^{12}})^n \frac{(1/2)_n(1/3)_n(2/3)_n(1/6)_n(5/6)_n}{(1)_n^5}.$$
	Let $$F(n,k) = \frac{2^{-2 k-6 n} \Gamma (-n-1) \Gamma \left(n-\frac{1}{2}\right) \Gamma (k+3 n+2) \Gamma \left(2 k+3 n+\frac{3}{2}\right)}{\Gamma \left(-n-\frac{5}{2}\right) \Gamma \left(-n-\frac{1}{2}\right) \Gamma (n) \Gamma (n+1) \Gamma (k+2 n+1) \Gamma (k+2 n+2) \Gamma (k+2 n+3)},$$
	conditions of Proposition \ref{WZ_prop} is satisfied, $g(n)=0$ and $\lim_{n\to\infty} \sum_{k\geq 0} F(n,k) = 0$, so we have $\sum_{k\geq 0} F(0,k+a) = \sum_{n\geq 1} G(n-1,a)$, that is
	\begin{multline*}\sum_{k\geq 0}\frac{2^{-6 a-2 k-1} (3 a+k+1) (6 a+4 k+1) (3 a+\frac{1}{2})_{2 k} (3 a+1)_k}{(2 a+k+1)^2 (2 a+k+2) ((2 a+1)_k)^3} \\
		+ \sum_{n\geq 0} \frac{(a+1) (2 a-1) (-1)^n 2^{-3 a-3 n-3/2} (2 a+2 n+1)^2 (2 a+2 n+3) (2 a+2 n+5) (6 a+6 n+7) \Gamma (2 a+1)^3 (a+\frac{1}{2})_n^3}{\sqrt{\pi} a (2 a+1)^2 (2 a+3) (2 a+5) (a+n+1)^2 (a+n+2) (2 a+2 n-1) \Gamma(3 a+\frac{1}{2}) \Gamma (3 a+1) (a+1)_n^3} \\
		= \sum_{n\geq 0} \frac{(a+1) (2 a-1) (-1)^n 2^{-6 a-6 n-11} (2 a+2 n+5) (3 a+3 n+1) (6 a+6 n+1)  (a+\frac{1}{2})_n^3 (3 a+\frac{1}{2})_{3 n} (3 a+1)_{3 n}\times P}{a (2 a+1)^2 (2 a+3) (2 a+5) (a+n+1)^3 (a+n+2) (2 a+2 n-1) (a+1)_n^3 (2 a+1)_{2 n}^3}.
	\end{multline*}
here {\small $P = 7720 a^4+30880 a^3 n+29988 a^3+46320 a^2 n^2+89964 a^2 n+41318 a^2+30880 a n^3+89964 a n^2+82636 a n+23499 a+7720 n^4+29988 n^3+41318 n^2+23499 n+4410$}. Comparing coefficient of $a^{-1}$ gives
\begin{multline*}-\sum_{n\geq 0} \frac{(-1)^n 2^{-3 n-\frac{3}{2}} (2 n+1)^2 (2 n+3) (2 n+5) (6 n+7) (\frac{1}{2})_n^3}{15 \pi  (n+1)^2 (n+2) (2 n-1) (1)_n^3} \\ = \sum_{n\geq 0} -\frac{(-1)^n 2^{-6 n-11} (2 n+5) (3 n+1) (6 n+1) (7720 n^4+29988 n^3+41318 n^2+23499 n+4410) (\frac{1}{2})_n^3 (\frac{1}{2})_{3 n} (1)_{3 n}}{15 (n+1)^3 (n+2) (2 n-1) (1)_n^3 (1)_{2 n}^3}.
\end{multline*}
The summand of RHS equals
$$\Delta_n\left(\frac{(-1)^n 2^{-6 n} n^2 (124 n^2+152 n+13)}{45(1+n)(2n-1)} a_n\right) + \frac{2^{-5} (45+549n+1930n^2)}{45}a_n;$$
the summand of LHS equals
$$\Delta_n\left(\frac{(-1)^n 2^{\frac{9}{2}-3 n} n (12 n^2+12 n-1) (\frac{1}{2})_n^3}{15 \pi  (n+1) (2 n-1) (1)_n^3}\right)+\frac{(-1)^n 2^{\frac{9}{2}-3 n} (6 n+1) (\frac{1}{2})_n^3}{15 \pi  (1)_n^3}.$$
A well-known result $\sum_{n\geq 0}\frac{(-1/8)^n (6 n+1) (\frac{1}{2})_n^3}{(1)_n^3} = \frac{2\sqrt{2}}{\pi}$ implies the conjecture. 
\end{romexample}

\begin{romexample}\label{12pi-2_ex}Z.W. Sun (\cite{sun2023new}) proposes that, $$\sum_{n\geq 0} a_n \frac{92n^3+54n^2+12n+1}{6n+1} = \frac{12}{\pi^2}, \qquad a_n = (\frac{4}{27})^n \frac{(1/2)_n^7}{(1/6)_n (5/6)_n (1)^5}.$$
	Consider $$F(n,k) = \frac{(4 k+4 n+3) \Gamma \left(n-\frac{1}{2}\right) \Gamma \left(k+n+\frac{1}{2}\right)^2 \Gamma \left(k+n+\frac{3}{2}\right) \Gamma \left(2 k+n+\frac{3}{2}\right) \Gamma (k+2 n+1)}{\Gamma \left(k+\frac{3}{2}\right) \Gamma \left(-n-\frac{1}{2}\right)^2 \Gamma (n)^2 \Gamma (n+1) \Gamma (k+n+1) \Gamma (k+n+2)^2 \Gamma \left(2 k+3 n+\frac{5}{2}\right)},$$
	one computes $G(n,k), g(n)\neq 0$ and $\lim_{n\to\infty} \sum_{k\geq 0} F(n,k) = 0$, $\sum_{k\geq 0} F(a,k) + \sum_{n\geq 0} g(n+a) = \sum_{n\geq 0} G(n+a,0)$ becomes
	\begin{multline*}\sum_{k\geq 0}\frac{(2 a+2 k+1) (4 a+4 k+3) (a+\frac{1}{2})_k^3 (a+\frac{1}{2})_{2 k+1} (2 a+1)_k}{(2 k+1) (a+k+1)^2 (\frac{1}{2})_k ((a+1)_k)^3 (3 a+\frac{1}{2})_{2 k+2}} \\ - \sum_{n\geq 0} \frac{3 \sqrt{\pi } (2 a-1) (-1)^{2 n} 2^{-2 a-2 n} (2 a+2 n+1)^2 \Gamma (a+1)^3 \Gamma (3 a+\frac{1}{2}) (a+\frac{1}{2})_n^3}{a^2 (2 a+1)^2 (a+n+1) (2 a+2 n-1) \Gamma (a+\frac{1}{2})^4 \Gamma (2 a+1) (a+1)_n^3} \\ = -\sum_{n\geq 0} \frac{\splitfrac{(2 a-1) (2 a+2 n+1)^3 (184 a^4+736 a^3 n+508 a^3+1104 a^2 n^2+1524 a^2 n+492 a^2}{+736 a n^3+1524 a n^2+984 a n+198 a+184 n^4+508 n^3+492 n^2+198 n+27) (a+\frac{1}{2})_n^7 (2 a+1)_{2 n}}}{6 a^2 (2 a+1)^2 (a+n+1)^3 (2 a+2 n-1) (6 a+6 n+1) (6 a+6 n+5) (a+1)_n)^6 (3 a+\frac{1}{2})_{3 n}}.
	\end{multline*}
	Extract coefficient of $a^{-2}$:
	$$\sum_{n\geq 0}\frac{3 (\frac{1}{2})^{2 n} (2 n+1)^2 (\frac{1}{2})_n^3}{\pi  (n+1) (2 n-1) (1)_n^3} = \sum_{n\geq 0}\frac{(-1)^{2 n} (2 n+1)^3 (184 n^4+508 n^3+492 n^2+198 n+27) (\frac{1}{2})_n^7 (1)_{2 n}}{6 (n+1)^3 (2 n-1) (6 n+1) (6 n+5) (\frac{1}{2})_{3 n} (1)_n^6}.
	$$
	The summand of RHS equals $$\Delta_n\left(\frac{4n^2 (1+6n)}{2n-1} a_n\right) + \frac{2(92n^3+54n^2+12n+1)}{3(6n+1)}a_n;$$
	while summand of LHS equals $$\Delta_n\left(\frac{2^{5-2 n} n^2 (\frac{1}{2})_n^3}{\pi  (2 n-1) (1)_n^3}\right)+\frac{2^{1-2 n} (6 n+1) (\frac{1}{2})_n^3}{\pi (1)_n^3}.$$
	A well-known result $\sum_{n\geq 0}\frac{2^{-2 n} (6 n+1) (\frac{1}{2})_n^3}{(1)_n^3} = \frac{4}{\pi}$ implies the conjecture. 
\end{romexample}

\begin{romexample}\label{896-zeta3_ex}We prove following conjecture appeared in \cite{guillera2016bilateral}, 
	$$\sum_{n\geq 1} (-\frac{2^8}{5^5})^n \frac{(1)_n^5}{(\frac{1}{2})_n(\frac{1}{5})_n(\frac{2}{5})_n(\frac{3}{5})_n(\frac{4}{5})_n} \frac{-483n^2+245n-30}{n^5} = 896\zeta(3).$$
	Consider 
	$$F(n,k) = \frac{(-4)^n (4 k+6 n+3) \Gamma (n+1)^2 \Gamma (2 n+1) \Gamma (2 k+n+1) \Gamma (k+2 n+1)^2}{\Gamma \left(n+\frac{1}{2}\right)^2 \Gamma \left(k+n+\frac{3}{2}\right)^2 \Gamma (2 k+5 n+3)},$$
	one find the corresponding $G(n,k)$, conditions of Proposition \ref{WZ_prop} is satisfied, $g(n)=0$ and $\lim_{n\to\infty} \sum_{k\geq 0} F(n,k) = 0$, so we have $\sum_{k\geq 0} F(0,k+a) = \sum_{n\geq 1} G(n-1,a)$, that is
	\begin{multline*}\sum_{k\geq 0}\frac{4 (6 a+4 k+3) (a+1)_{2 k} (2a+1)_k^2}{(2 a+2 k+1)^2 (a+\frac{1}{2})_k^2 (5 a+1)_{2 k+2}} = \\ 
		\sum_{n\geq 1} \frac{(-1)^{n+1} 2^{2 n-7} (483 a^2+966 a n-245 a+483 n^2-245 n+30)(a+1)_n^3 (2 a+1)_{2 n}^3}{(a+n)^5 (a+1/2)_n^4 (5 a+1)_{5 n}}.
	\end{multline*}
	When $a=0$, RHS recovers exactly the sum concerned, and it is easy to show $\sum_{k\geq 0} \frac{4k+3}{(2k+1)^3 (k+1)} \frac{(1)_k^2}{(1/2)_k^2} = \frac{7}{2}\zeta(3)$, this proves the formula. 
\end{romexample}

\begin{romexample}\label{ex_zeta5}Following marvellous series of $\zeta(5)$ was empirically discovered by Y.Zhao\footnote{\url{https://mathoverflow.net/questions/281009}}.
	$$\sum_{n\geq 1}\left(-\frac{2^{10}}{5^5}\right)^n \frac{(1)_n^9}{(\frac{1}{2})_n^5(\frac{1}{5})_n(\frac{2}{5})_n(\frac{3}{5})_n(\frac{4}{5})_n} \frac{5532n^4-5600n^3+2275n^2-425n+30}{n^9} =-380928\zeta(5).$$
	
	Let $$F(n,k) = \frac{(-1)^n (4 k+6 n+3) \Gamma (n+1)^4 \Gamma (2 n+1) \Gamma (2 k+n+1) \Gamma (k+2 n+1)^4}{\Gamma \left(n+\frac{1}{2}\right)^6 \Gamma \left(k+n+\frac{3}{2}\right)^4 \Gamma (2 k+5 n+3)},$$
	apply Proposition \ref{WZ_prop} (here $g(n)\neq 0$), $\sum_{k\geq 0} F(0,k) + \sum_{n\geq 0} g(n) = \sum_{n\geq 1} G(n-1,0)$ becomes
\begin{multline*}\sum_{k\geq 0} \frac{8 (4 k+3)(1)_k^4}{(k+1) (2 k+1)^5 (1/2)_k^4} + \pi ^2\sum_{n\geq 0} \frac{ (-1)^n 4^{1-2 n} (10 n^2+14 n+5) (1)_{2n} (1)_n^4}{(2 n+1)^5 (1/2)_n^6} \\ = -\sum_{n\geq 1} \frac{(-1)^n(5532 n^4-5600 n^3+2275 n^2-425 n+30) \Gamma (n)^5 \Gamma (2 n)^5}{320 (1/2)_n^{10} \Gamma (5 n)}.\end{multline*}

We will show in Example \ref{ex_reg} below, the summation in $k$ equals $2( 93 \zeta (5)-7 \pi ^2 \zeta(3))$, this is not trivial. The sum in $n$ at LHS (coming from $\sum_{n\geq 0}g(n)$) equals $14 \zeta(3)$, this is a known WZ-type conclusion. Combining these two evaluations give the result. 
\end{romexample}

\begin{remark}
This formula of $\zeta(5)$ has been used by D. Broadhurst (\cite{broadhurst2024five}) to calculate its numerical value to 200 billion decimal digits. 
\end{remark}

\begin{remark}This formula for $\zeta(5)$ seems to be the only currently known (conjectural or proven) formula for $\zeta(5) = \sum_{n\geq 1} u_n$ with
\begin{itemize}
\item $u_n$ is hypergeometric (i.e. $u_{n+1}/u_n \in \overline{\mathbb{Q}}(n)$), and
\item $u_n$ decreases at geometric rate (so trivial examples $\sum_{n\geq 1} (\pm 1)^n n^{-5}$ are excluded)
\end{itemize}
For $\zeta(3)$, examples of such representations are long known (e.g. Apéry's series). 
\end{remark}

\begin{romexample}\label{ex_pi-3}
	Let us establish following two sums due to B. Gourevich and J. Guillera (\cite{zudilin2011arithmetic}, \cite{guillera2008hypergeometric}, \cite{cohen2021rational}): \begin{align}\label{ex_pi-3_eq2}&\sum_{n\geq 0} \frac{(1/2)_n^7}{(1)_n^7} \frac{1}{2^{6n}} (168n^3+76n^2+14n+1) = \frac{32}{\pi^3}, \\
		&\sum_{n\geq 1} \frac{(1)_n^7}{(1/2)_n^7} \frac{1}{2^{6n}} \frac{21 n^3-22 n^2+8 n-1}{n^7} = \frac{\pi^4}{8}. \end{align}
	
We take $$F_1(n.k) = \frac{(-1)^n (4 k+6 n+3) \Gamma \left(n+\frac{1}{2}\right)^6 \Gamma (n+1)^4 \Gamma \left(k+n+\frac{1}{2}\right) \Gamma (k+n+1)^4 \Gamma (k+4 n+1)}{\Gamma (2 n) \Gamma \left(2 n+\frac{1}{2}\right)^4 \Gamma \left(k-n+\frac{3}{2}\right) \Gamma \left(k+2 n+\frac{3}{2}\right)^4 \Gamma (k+2 n+2)}.$$
	
	Computing required data $G_1, g_1$ of Proposition \ref{WZ_prop}, we obtain (here $g_1(n)\neq 0$ and $\lim_{n\to\infty} \sum_{k\geq 0} F_1(n,k) = 0$), 
	\begin{multline}\label{ex_pi-3_eq1} \sum_{k\geq 0} \frac{-32(6 a+4 k+3)(a+\frac{1}{2})_k (a+1)_k^4 (4 a+1)_k}{(2 a-2 k-1) (2 a+k+1) (4 a+2 k+1)^4 (\frac{1}{2}-a)_k (2 a+\frac{1}{2})_k^4 (2 a+1)_k} \\ - \frac{16 \Gamma \left(\frac{1}{2}-a\right) \Gamma \left(2 a+\frac{1}{2}\right)^4 \Gamma (2 a+1)}{a \Gamma \left(a+\frac{1}{2}\right) \Gamma (a+1)^4 \Gamma (4 a+1)} \sum_{n\geq 0} \frac{(-1)^n p_1(a,n) (a+\frac{1}{2})_n^6 (a+1)_n^4}{(4 a+4 n+1)^4 (4 a+4 n+3)^4 (2 a+\frac{1}{2})_{2 n}^4 (2 a+1)_{2 n}} \\ = -\sum_{n\geq 1}\frac{p_2(a,n) (a+\frac{1}{2})_n^8 (a+1)_n^8 (4 a+1)_{4 n}}{16 a (a+n)^7 (2 a+2 n-1)^7 (2 a+\frac{1}{2})_{2 n}^8 (2 a+1)_{2 n}^2}\end{multline}
	here {\small $p_1(a,n) = 1640 a^6+9840 a^5 n+5936 a^5+24600 a^4 n^2+29680 a^4 n+8738 a^4+32800 a^3 n^3+59360 a^3 n^2+34952 a^3 n+6664 a^3+24600 a^2 n^4+59360 a^2 n^3+52428 a^2 n^2+19992 a^2 n+2762 a^2+9840 a n^5+29680 a n^4+34952 a n^3+19992 a n^2+5524 a n+587 a+1640 n^6+5936 n^5+8738 n^4+6664 n^3+2762 n^2+587 n+50$} and $p_2(a,n) \in \mathbb{Z}[a,n]$ is even longer\footnote{readers interested in the full expression should perform Gosper's algorithm themselves.}. With the full expression of $p_2(a,n)$, one easily verifies, by $\sum_{n\geq 1} a_n = \sum_{n\geq 1} a_{2n}+a_{2n-1}$, the sum on RHS actually bisects, and equals
	$$-\sum_{n\geq 1} \frac{2^{3-8 n} (168 a^3+252 a^2 n-88 a^2+126 a n^2-88 a n+16 a+21 n^3-22 n^2+8 n-1) (2 a+1)_n^6 (4 a+1)_{2 n}}{a (2 a+n)^7 (2 a+\frac{1}{2})_n^8}.$$
	
	Now we shall see that $\sum_{n\geq 0} g_1(n)$ from above is actually $\sum_{n\geq 0} G_2(n,0)$ of another WZ-pair\footnote{actually this is not a coincidence, the $g(n)$ from balanced WZ-seed \textsf{Dougall7F6} is always the $G(n,k)$ of another WZ-pair coming from seed \textsf{Dougall5F4}.} $(F_2, G_2)$:
	$$F_2(n,k) = \frac{(3 a+2 k+3 n+1) \Gamma \left(a+n+\frac{1}{2}\right)^2 \Gamma (a+n+1)^3 \Gamma \left(a+k+n+\frac{1}{2}\right)^3 \Gamma (a+k+n+1)}{\Gamma \left(-a-n+\frac{1}{2}\right) \Gamma \left(2 a+2 n+\frac{1}{2}\right) \Gamma (2 a+k+2 n+1) \Gamma \left(2 a+k+2 n+\frac{3}{2}\right)^3}.$$
	then corresponding $g_2(n) = 0$ and $\lim_{n\to\infty}\sum_{k\geq 0} F_2(n,k) = 0$, $$\sum_{k\geq 0} F_2(0,k) = \sum_{n\geq 0} G_2(n,0),$$ becomes
	$$\sum_{k\geq 0} \frac{(3 a+2 k+1) (a+\frac{1}{2})_k^3 (a+1)_k}{2(4 a+2 k+1)^3 (2 a+\frac{1}{2})_k^3 (2 a+1)_k} = \sum_{n\geq 0} \frac{(-1)^n p_1(a,n) (a+\frac{1}{2})_n^6 (a+1)_n^4}{(4 a+4 n+1)^4 (4 a+4 n+3)^4 (2 a+\frac{1}{2})_{2 n}^4 (2 a+1)_{2 n}}.$$
	
	Combining above observations, equation (\ref{ex_pi-3_eq1}) becomes
	\begin{multline}\label{ex_pi-3_eq3}\sum_{k\geq 0} \frac{-32(6 a+4 k+3)(a+\frac{1}{2})_k (a+1)_k^4 (4 a+1)_k}{(2 a-2 k-1) (2 a+k+1) (4 a+2 k+1)^4 (\frac{1}{2}-a)_k (2 a+\frac{1}{2})_k^4 (2 a+1)_k} \\ - \frac{8 \Gamma \left(\frac{1}{2}-a\right) \Gamma \left(2 a+\frac{1}{2}\right)^4 \Gamma (2 a+1)}{a \Gamma \left(a+\frac{1}{2}\right) \Gamma (a+1)^4 \Gamma (4 a+1)} \sum_{k\geq 0} \frac{(3 a+2 k+1) (a+\frac{1}{2})_k^3 (a+1)_k}{(4 a+2 k+1)^3 (2 a+\frac{1}{2})_k^3 (2 a+1)_k} \\ = -\sum_{n\geq 1} \frac{2^{3-8 n} (168 a^3+252 a^2 n-88 a^2+126 a n^2-88 a n+16 a+21 n^3-22 n^2+8 n-1) (2 a+1)_n^6 (4 a+1)_{2 n}}{a (2 a+n)^7 (2 a+\frac{1}{2})_n^8}.\end{multline}
	Comparing coefficient of $a^{-1}$ to which first term on LHS does not contribute:
	$$-8\pi^2\sum_{k\geq 0} \frac{1}{(2k+1)^2} = -\sum_{n\geq 1}\frac{2^{3-8 n} (21 n^3-22 n^2+8 n-1) (1)_n^6 (1)_{2 n}}{n^7 (\frac{1}{2})_n^8};$$
	This is equivalent to equation (\ref{ex_pi-3_eq2}). Now consider Laurent expansion of equation (\ref{ex_pi-3_eq3}) at $a=-1/4$, and look at coefficient of leading term $(a+1/4)^{-7}$, again, the first term does not contribute and we get
	\begin{multline*}\lim_{a\to 1/4} \left((a+1/4)^7 \frac{ \Gamma \left(\frac{1}{2}-a\right) \Gamma \left(2 a+\frac{1}{2}\right)^4 \Gamma (2 a+1)}{a \Gamma \left(a+\frac{1}{2}\right) \Gamma (a+1)^4 \Gamma (4 a+1)} \sum_{k\geq 0} \frac{(3 a+2 k+1) (a+\frac{1}{2})_k^3 (a+1)_k}{(4 a+2 k+1)^3 (2 a+\frac{1}{2})_k^3 (2 a+1)_k} \right) \\ = -\sum_{n\geq 0} \frac{2^{-8 n} (6 n-5) (28 n^2-48 n+21) (\frac{1}{2})_n^6 (1)_{2 n-1}}{(2 n-1)^7 (1)_{n-1}^8}.\end{multline*}
	It is easy to show LHS equals $-2^{-9}\pi^{-3}$, while RHS is equivalent to the sum expression in equation (\ref{ex_pi-3_eq2}).
\end{romexample}

\begin{remark}Another proof of above $1/\pi^3$-formula using a simpler WZ-pair can be found in appendix of the recently published article \cite{au2024wilfQ}.\end{remark}

\begin{romexample}\label{-2beta4_ex}
	We establish, here $L_{-4}(s)$ is Dirichlet's beta function $$\begin{aligned}\sum_{n\geq 1} (\frac{-1}{8})^n \frac{(1)_n^3}{(1/2)_n^3} \frac{3n-1}{n^3} &= -2 L_{-4}(2), \\
		\sum_{n\geq 1} (\frac{-1}{8})^n \frac{(1)_n^3}{(1/2)_n^3} \frac{2 (3 n-1) (H_{2 n-1}-H_{n-1})-1}{n^3} &= 4 \Im\left(\text{Li}_3\left(\frac{1}{2}+\frac{i}{2}\right)\right)-\frac{5 \pi ^3}{32}-\frac{1}{8} \pi  \log ^2(2), \\
		\sum_{n\geq 1} (\frac{-1}{8})^n \frac{(1)_n^3}{(1/2)_n^3} \frac{(3 n-1) (H_{2 n-1}^{(2)}-\frac{5 H_{n-1}^{(2)}}{4})}{n^3} &= -2 L_{-4}(4), \\
		\sum_{n\geq 1} (\frac{-1}{8})^n \frac{(1)_n^3}{(1/2)_n^3} \frac{(3 n-1) (\frac{7 H_{n-1}^{(3)}}{8}+H_{2 n-1}^{(3)})}{n^3} &= 3\Im \text{Li}_{4,1}(i,1)+\Im \text{Li}_{4,1}(i,-1)+2 L_{-4}(2) \log (2)-\frac{5 \pi ^5}{512}.\end{aligned}$$
	The first was proved by Guillera in \cite{guillera2008hypergeometric} and third equality was conjectured by Z.W. Sun\footnote{He also conjectured $p$-adic analogues of second and fourth equalities, but was unable to guess their results in the Archimedean case} in \cite{sun2022conjectures} and proved by Wei in \cite{wei2023conjectural}. RHS falls out naturally once we show certain infinite series gives level 4 CMZVs when extracting coefficients. 
	
	Let $$F(n,k) = \frac{2^{-2 k-3 n} \Gamma (-a+2 b+n+1) \Gamma (a+2 k+n+1)}{\Gamma (a-n) \Gamma \left(k+n+\frac{3}{2}\right) \Gamma \left(-a+b+n+\frac{1}{2}\right) \Gamma \left(b+k+n+\frac{3}{2}\right)},$$ then in Proposition \ref{WZ_prop}, $\lim_{n\to\infty} \sum_{k\geq 0} F(n,k) = 0$ and $g(n)=0$, $\sum_{k\geq 0} F(0,k+c) = \sum_{n\geq 1} G(n-1,c)$ becomes
	\begin{multline*}\frac{1}{(1+2c)(1+2b+2c)} \pFq{3}{2}{1,\frac{a}{2}+c+\frac{1}{2},\frac{a}{2}+c+1}{c+\frac{3}{2},b+c+\frac{3}{2}}{1} \\ = \sum_{n\geq 1}\frac{(-1)^{n+1} 2^{-3 n} (1-a)_n (-a+2 b+2 c+3 n-1) (-a+2 b+1)_n (a+2 c+1)_n}{(n-a) (-a+2 b+n) (a+2 c+n) (c+\frac{1}{2})_n (-a+b+\frac{1}{2})_n (b+c+\frac{1}{2})_n}.\end{multline*}
	where we rewrote the sum over $k$ using hypergeometric function, it's not obvious whether coefficient of $a^i b^j c^k$ are CMZVs of some level since condition (\ref{int_diff_cond}) is not met. Similar to Example \ref{ex_CMZVtrans}, we try to apply some hypergeometric transformation to it:  Whipple's transformation 
	$$\pFq{6}{5}{a,\frac{a}{2}+1,b,c,d,e}{\frac{a}{2},a-b+1,a-c+1,a-d+1,a-e+1}{-1} =\frac{\Gamma (a-d+1) \Gamma (a-e+1)}{\Gamma (a+1) \Gamma (a-d-e+1)} \pFq{3}{2}{a-b-c+1,d,e}{a-b+1,a-c+1}{1},$$
	transforms above $_3F_2$ to
	$$ \frac{\Gamma(-a+b+\frac{1}{2}) \Gamma (b+2 c+2)}{\Gamma \left(-\frac{a}{2}+b+c+1\right) \Gamma \left(-\frac{a}{2}+b+c+\frac{3}{2}\right)}\pFq{6}{5}{b+2 c+1,\frac{b}{2}+c+\frac{3}{2},b+c+\frac{1}{2},c+\frac{1}{2},\frac{a}{2}+c+\frac{1}{2},\frac{a}{2}+c+1}{\frac{b}{2}+c+\frac{1}{2},c+\frac{3}{2},b+c+\frac{3}{2},-\frac{a}{2}+b+c+\frac{3}{2},-\frac{a}{2}+b+c+1}{-1}.$$
	In this form, equation (\ref{int_diff_cond}) holds and we conclude coefficients of $a^i b^j c^k$ are CMZVs of weight $a+b+c+2$ and level $4$. 
\end{romexample}

	\section{Examples with non-vanishing boundary term}
	We give some examples such that $$\lim_{n\to\infty} \sum_{k\geq 0}F(n,k) \neq 0$$ in Proposition \ref{WZ_prop}. Such limits, whose existence and finiteness are granted by the Proposition's statement, usually cannot be calculated naively using dominated convergence or anything that tries to move limit inside summation. Efforts on asymptotic analysis seem sometimes necessary, but techniques we shall employ are all standard (e.g. can be found in \cite{de1981asymptotic}). 
	
	Two such examples are already in the author's previous paper: Example IX and X of \cite{au2022multiple}. We give some more examples below.
\begin{romexample}\label{ex_nonvanishing_1}
		Z.W. Sun conjectured\footnote{\url{https://mathoverflow.net/questions/457326}} that $$\sum_{n\geq 0} (\frac{27}{32})^n \frac{(1/6)_n (5/6)_n}{(1)_n^2} \frac{10n-1}{2n+1} = 0.$$
		
	By searching the exponent $27/32$ and hypergeometric part $\frac{(1/6)_n (5/6)_n}{(1)_n^2}$, we are inspired to take WZ-candidate
 	$$F(n,k)=  \frac{2^{-2 k-3 n} \Gamma \left(-a+c-n+\frac{1}{2}\right) \Gamma \left(a+c+2 d+2 k+3 n+\frac{1}{2}\right)}{\Gamma \left(a+c-n+\frac{1}{2}\right) \Gamma (d+k+2 n+1) \Gamma (c+d+k+n+1)}.$$
	One its WZ-mate $G$, and checks prerequisites of Proposition \ref{WZ_prop} are satisfied (with $g(n)=0$) when $a,c,d$ are in a neighbourhood of $0$ with $\Re(a)<0$. We have $$\mathcal{E}_F(x) = 2^{-2 x-3} | x+1| ^{-x-1} | x+2| ^{-x-2} | 2 x+3| ^{2 x+3},$$
	and $\lim_{x\to \infty} \mathcal{E}_F(x) = 1$, so Proposition \ref{vanishing_prop} is not satisfied. It indeed happens that the limit is not zero. 
	
	After some rearrangements, one has\begin{equation}\label{eq_1}\sum_{k\geq 0} \frac{2^{-2 k} \left(a+c+2 d+\frac{1}{2}\right)_{2 k}}{(d+1)_k (c+d+1)_k} = L -\frac{1}{4}\times \sum_{n\geq 0} \frac{2^{-3 n} P\times (-a-c+\frac{1}{2})_n (a+c+2 d+\frac{1}{2})_{3 n}}{(d+2 n+1) (2 a-2 c+2 n+1) (d+1)_{2 n} \left(a-c+\frac{1}{2}\right)_n (c+d+1)_n}.\end{equation}
	here $\small{P = 4 a^2+8 a c+8 a d+8 a n+4 c^2+24 c d+40 c n+16 c+16 d^2+40 d n+12 d+20 n^2+8 n-1}$ and $$L = \frac{\Gamma(a+c+1/2)\Gamma(d+1)\Gamma(c+d+1)}{\Gamma(a+c+2d+1/2) \Gamma(1/2-a+c)} \times \lim_{n\to\infty} \sum_{k\geq 0} \frac{2^{-2 k-3 n}\Gamma(-a+c-n+\frac{1}{2}) \Gamma(a+c+2 d+2 k+3 n+\frac{1}{2})}{\Gamma \left(-a+c+\frac{1}{2}\right)\Gamma \left(a+c-n+\frac{1}{2}\right) \Gamma (d+k+2 n+1) \Gamma (c+d+k+n+1)}.$$
	We know $L$ exists and is finite when $a,c,d$ are near $0$ and $\Re(a)<0$, we show below $$L = 2^{-a+c+2 d-1/2} \frac{\Gamma (-a) \Gamma (d+1)  \Gamma \left(a-c+\frac{1}{2}\right) \Gamma (c+d+1)}{\sqrt{\pi } \Gamma \left(-a-c+\frac{1}{2}\right) \Gamma \left(a+c+2 d+\frac{1}{2}\right)}.$$
	
	When $d=-c$, LHS of \ref{eq_1} can be summed using Gauss $_2F_1$ theorem, it can then be seen to equal to $L$. Therefore for $a,c\in \mathbb{C}$, 
	$$\sum_{n\geq 0}\frac{2^{-3 n} \left(4 a^2+8 a n-4 c^2+4 c+20 n^2+8 n-1\right) \left(-a-c+\frac{1}{2}\right)_n \left(a-c+\frac{1}{2}\right)_{3 n}}{(-c+2 n+1)(2 a-2 c+2 n+1)  (1)_n (1-c)_{2 n} \left(a-c+\frac{1}{2}\right)_n} = 0,$$
	
	one specializes to $a=c=0$ to get Sun's conjecture.

\begin{lemma}
The limit $L$ in above example equals $$2^{-a+c+2 d-1/2} \frac{\Gamma (-a) \Gamma (d+1)  \Gamma \left(a-c+\frac{1}{2}\right) \Gamma (c+d+1)}{\sqrt{\pi } \Gamma \left(-a-c+\frac{1}{2}\right) \Gamma \left(a+c+2 d+\frac{1}{2}\right)}.$$
\end{lemma}
\begin{proof}[Proof sketch]
$L$ can be summed in terms of $_3F_2$ with results in form
$$L = (\text{some gamma factors in }n) \times \pFq{3}{2}{1,1/4+a/2+c/2+d+3n/2,3/4+a/2+c/2+d+3n/2}{1 + c + d + n, 1 + d + 2 n}{1}$$
the asymptotic of gamma product is trivial to find via Stirling's formula, so we only concentrate on the term whose leading asymptotic is non-trivial to find, i.e. the $_3F_2$ term. We first transform the $_3F_2$ into integral, it equals
\begin{multline*}\frac{\Gamma (d+2 n+1) \Gamma (c+d+n+1)}{\Gamma \left(-\frac{a}{2}+\frac{n}{2}-\frac{c}{2}+\frac{1}{4}\right) \Gamma (c+d+n) \Gamma \left(\frac{a}{2}+\frac{c}{2}+d+\frac{3 n}{2}+\frac{3}{4}\right)} \times \\ \int_{[0,1]^2} (1-t_1)^{-3/4-a/2-c/2} t_1^{-1/4+a/2+c/2+d} (1-t_2)^{-1+c+d} (1-t_1t_2)^{-1/4-a/2-c/2-d} \left( \frac{(1-t_1)t_1^3 (1-t_2)}{(1-t_1t_2)^3}\right)^{n/2} dt_1 dt_2.\end{multline*}
Denote the integral by $I$, the main contribution of $I$ comes from the neighbourhood of $(t_1,t_2) = (1,1)$. We first make some substitutions so that this maximum will be attained in interior of integration region, to be precise, successively let $t_i = 1/(1+s_i)$, $s_2 \mapsto s_2\times s_1$, the main contribution will now be around $s_1 = 0$, and $I$ becomes
$$\int_{[0,\infty)^2} s_1^{-1-a} (1+s_1)^{-3/4+a/2+c/2} s_2^{-1+c+d} (1+s_1s_2)^{-3/4+a/2-c/2} (1+s_2 + s_1s_2)^{-1/4-a/2-c/2-d} \left(\frac{s_2^2 (1+s_1s_2)}{(1+s_2(1+s_1))^3}\right)^{n/2} ds_i.$$
Recall we can assume $s_1$ to be near $0$, when this is the case, the function (with $s_1$ fixed) $\frac{s_2^2 (1+s_1s_2)}{(1+s2(1+s_1))^3}$ has a unique maximum at $s_2 = 2/(1-2s_1) := s_0 \in (0,\infty)$, so we can apply the Laplace method $$\int_0^\infty h(s_2) \exp(n g(s_2)) ds_2 \sim \sqrt{\frac{2\pi}{n |g''(s_0)|}} h(s_0)\exp(n g(s_0))$$ to get the leading term of integral w.r.t. to $s_2$. After some calculation, we have
$$I\sim \int_0^\infty \frac{\sqrt{\pi } s_1^{-a-1} 2^{c+d+n+\frac{1}{2}} 3^{\frac{1}{4} (-2 a-2 c-4 d-6 n+1)} (s_1+1)^{\frac{1}{4} (2 a+2 c-2 n-3)}}{\sqrt{n}} ds_1.$$
Now finding leading term is easy since the integral is beta function. Plug it back and calculate asymptotic of all gamma products gives the result.
\end{proof}
\end{romexample}

\begin{romexample}\label{ex_nonvanishing_2}
	Z.W. Sun in \cite{sun2023new} conjectures the following\footnote{first equality was already established by the author using a two-variable hypergeometric identity (\cite{au2022multiple}, Example X), which is unfortunately not strong enough to derive the second identity.}: 
	$$\sum_{n\geq 1} (\frac{3}{4})^n \frac{(1)_n^2}{(1/4)_n (3/4)_n} \frac{1}{n(2n-1)} = \frac{15}{2}L_{-3}(2),$$
	$$\sum_{n\geq 1} (\frac{3}{4})^n \frac{(1)_n^2}{(1/4)_n (3/4)_n} \frac{6H_{4n-1}-9H_{2n-1}+2H_{n-1}+6/(2n-1)}{n(2n-1)} = \frac{2\pi^3}{\sqrt{3}}.$$
	Actually, both follow from a three-variable identity:
	\begin{multline*}\sum_{k\geq 0} \frac{(1+a+2c+2k) (1+a+c)_k (1/2+b+3c)_{3k}}{(1+c)_k (1/2+3a-b+3c)_{3k+3}} = \frac{3^{-b-\frac{3}{2}} \Gamma \left(b+\frac{1}{2}\right) \Gamma (c+1) \Gamma \left(a-b+\frac{1}{2}\right) \Gamma \left(3 a-b+3 c+\frac{1}{2}\right)}{\Gamma (3 a-2 b+2) \Gamma (a+c+1) \Gamma \left(b+3 c+\frac{1}{2}\right)} \\ - \frac{2c}{1+3a-2b} \sum_{n\geq 1} \frac{3^{n-1} (3 a-2 b+1)_n (a+c+1)_n \left(b+3 c+\frac{1}{2}\right)_n}{(3 a-2 b+1) (a+c+n) (2 b+6 c+2 n-1) (b+\frac{1}{2})_n (3 a-b+3 c+\frac{1}{2})_{2 n}}.\end{multline*}
	For example, the first conjecture follows by taking $a^0b^0c^1$ coefficient. Note that coefficients of $a^ib^jc^k$ of LHS are level 6 CMZVs. The equality itself is again a direct application of Proposition \ref{WZ_prop} with $$F(n,k) = \frac{3^{n} (a+2 k+n+1) \Gamma (3 a-2 b+n+2) \Gamma (a+k+n+1) \Gamma \left(b+3 k+n+\frac{1}{2}\right)}{\Gamma (k+1) \Gamma \left(a-b+\frac{1}{2}\right) \Gamma \left(b+n+\frac{1}{2}\right) \Gamma \left(3 a-b+3 k+2 n+\frac{7}{2}\right)},$$
	one verify $G(n,k)$ exists and $g(n)=0$, after some calculations, $$\sum_{k\geq 0} F(0,c+k) = \lim_{n\to\infty} \sum_{k\geq 0} F(n,k+c) + \sum_{n\geq 0} G(n,c)$$ becomes exactly above equality, except it's not obvious the limit equals the gamma product term on RHS, here we need to find
	\begin{multline*}\lim_{n\to\infty} \frac{3^n \Gamma \left(b+\frac{1}{2}\right) \Gamma (c+1) \Gamma \left(3 a-b+3 c+\frac{1}{2}\right) \Gamma (3 a-2 b+n+2) }{(3 a-2 b+1) \Gamma (3 a-2 b+1) \Gamma (a+c+1) \Gamma \left(b+3 c+\frac{1}{2}\right) \Gamma \left(b+n+\frac{1}{2}\right)  } \\ \times \sum_{k\geq 0} \frac{(a+2 c+2 k+n+1) \Gamma (a+c+k+n+1) \Gamma \left(b+3 c+3 k+n+\frac{1}{2}\right)}{\Gamma (c+k+1) \Gamma \left(3 a-b+3 c+3 k+2 n+\frac{7}{2}\right)}.\end{multline*} 
	
	\begin{lemma}\label{lemma_hypertrans}
	Above limit equals $\frac{3^{-b-\frac{3}{2}} \Gamma \left(b+\frac{1}{2}\right) \Gamma (c+1) \Gamma \left(a-b+\frac{1}{2}\right) \Gamma \left(3 a-b+3 c+\frac{1}{2}\right)}{\Gamma (3 a-2 b+2) \Gamma (a+c+1) \Gamma \left(b+3 c+\frac{1}{2}\right)}$.
	\end{lemma}
\begin{proof}[Proof sketch]
This can be shown using same method as above lemma, but we present an alternative method using transformation of very-well-poised $_7F_6$. We concentrate on the asymptotic of the hypergeometric sum, it is a shift (by $c$) of a very-well-poised $_5F_4$, so is a very-well-poised $_7F_6$, explicitly, it is
$$ \pFq{7}{6}{\frac{a}{2}+c+\frac{n}{2}+\frac{3}{2},1,a+c+n+1,\frac{b}{3}+c+\frac{n}{3}+\frac{1}{6},\frac{b}{3}+c+\frac{n}{3}+\frac{1}{2},\frac{b}{3}+c+\frac{n}{3}+\frac{5}{6}}{\frac{a}{2}+c+\frac{n}{2}+\frac{1}{2},c+1,a-\frac{b}{3}+c+\frac{2 n}{3}+\frac{11}{6},a-\frac{b}{3}+c+\frac{2 n}{3}+\frac{3}{2},a-\frac{b}{3}+c+\frac{2 n}{3}+\frac{7}{6}}{1}.$$
We hope, under transformation \ref{VWP7F6trans}, there will be a form such that the asymptotic with $n\to\infty$ is easy to see\footnote{to be precise, we want a form in which the number of numerator parameters of $_7F_6$ in which $n$ appear has to be the same as that of denominator parameters}. Indeed, use the formula, the above very-well-poised $_7F_6$ equals
\begin{multline*}(\text{some gamma factors}) \\ \times \pFq{7}{6}{a-b+c+\frac{1}{2},\frac{a}{2}-\frac{b}{2}+\frac{c}{2}+\frac{5}{4},c,-\frac{b}{3}-\frac{n}{3}+\frac{5}{6},-\frac{b}{3}-\frac{n}{3}+\frac{1}{2},-\frac{b}{3}-\frac{n}{3}+\frac{1}{6},a-b+\frac{1}{2}}{\frac{a}{2}-\frac{b}{2}+\frac{c}{2}+\frac{1}{4},a-b+\frac{3}{2},a-\frac{2 b}{3}+c+\frac{n}{3}+\frac{2}{3},a-\frac{2 b}{3}+c+\frac{n}{3}+1,a-\frac{2 b}{3}+c+\frac{n}{3}+\frac{4}{3},c+1}{1},\end{multline*}
as $n\to \infty$, above $_7F_6$ tends to
$$\pFq{4}{3}{a-b+c+\frac{1}{2},\frac{a}{2}-\frac{b}{2}+\frac{c}{2}+\frac{5}{4},c,a-b+\frac{1}{2}}{\frac{a}{2}-\frac{b}{2}+\frac{c}{2}+\frac{1}{4},a-b+\frac{3}{2},c+1}{-1},$$
which can be easily shown to be $\frac{\Gamma (c+1) \Gamma \left(a-b+\frac{3}{2}\right)}{\Gamma \left(a-b+c+\frac{3}{2}\right)}$. Go back and calculate asymptotic of all gamma products which we picked up along the way will give the result.
\end{proof}
\end{romexample}

\begin{romexample}\label{ex_reg}
Let $$F(n,k) = \frac{\splitfrac{\left(a+2 (k+n)+\frac{1}{2}\right) \Gamma (c+k) \Gamma (a-c-d+1) \Gamma \left(a+k+n+\frac{1}{2}\right) \Gamma \left(b+k+n+\frac{1}{2}\right)}{ \Gamma \left(d+k+2 n+\frac{1}{2}\right) \Gamma (a-b-c+n+1) \Gamma \left(a-b-d-n+\frac{1}{2}\right)}}{\splitfrac{\Gamma \left(b+\frac{1}{2}\right) \Gamma (c-n) \Gamma \left(d+n+\frac{1}{2}\right) \Gamma (k+n+1) \Gamma (a-d+k+1)}{ \Gamma \left(a-b-c-d+\frac{1}{2}\right) \Gamma (a-b+k+n+1) \Gamma \left(a-c+k+2 n+\frac{3}{2}\right)}},$$
this comes from the WZ-seed \textsf{Dougall5F4}. Assuming $a,b,c,d,e$ are in a neighbourhood of $0$, we are interested with the formula induced from this $F(n,k)$:
$$\sum_{k\geq 0} F(0,k+e+1/2) = \lim_{n\to \infty}\sum_{k\geq 0} F(n,k+e+1/2) + \sum_{n\geq 1} G(n-1,e+1/2).$$
One checks $\sum_{n\geq 1} G(n-1,e+1/2)$ converges if $\Re(a-c-d)<0$. However, we encounter three difficulties when writing out this formula explicitly: \par
First is non-vanishing of $\lim_{n\to \infty}\sum_{k\geq 0} F(n,k+e+1/2)$, it can be found using the same technique as in Lemma \ref{lemma_hypertrans}: the sum over $k$ is a very-well-poised $_7F_6$ and suitable transformation will mould it into a form in which the limit is easily seen. \par
Second is, when comparing coefficient of LHS at $(a,b,c,d,e)=(0,0,0,0,0)$, it's not obvious its coefficients are CMZVs, this can be resolved using the same technique as in Example (\ref{ex_CMZVtrans}): suitable transformation (\ref{VWP7F6trans}) will put it into a form in which condition (\ref{int_diff_cond}) holds. \par
 After doing above two steps, we obtain the following formula:
	\begin{equation}\label{ex_eq1}\sum_{k\geq 0} A(k) = B + \sum_{n\geq 1} C(n),\end{equation}
	with \begin{multline*}A(k) = \frac{\Gamma \left(e+\frac{1}{2}\right) \Gamma \left(a-b+e+\frac{1}{2}\right) \Gamma (a-c+e+1) \Gamma (2 a-b-c-2 d+e+1)}{\Gamma (d+e+1) \Gamma (a-b-d+e+1) \Gamma \left(a-c-d+e+\frac{1}{2}\right) \Gamma \left(2 a-b-c-d+e+\frac{1}{2}\right)} \\ \times \frac{4 (2 a-b-c-2 d+e+2 k+1) \left(\frac{1}{2}-d\right)_k \left(a-b-d+\frac{1}{2}\right)_k (a-c-d+1)_k  (2 a-b-c-2 d+e+1)_k}{\splitfrac{(1)_k (2 a-2 d+2 e+2 k+1) (2 a-2 b-2 c-2 d+2 k+1)}{ (a-b-d+e+1)_k \left(a-c-d+e+\frac{1}{2}\right)_{k+1} \left(2 a-b-c-d+e+\frac{1}{2}\right)_{k+1}}}, \end{multline*}
	here $\sum_{k\geq 0} A(k)$, when expanded around $(a,b,c,d,e)=(0,0,0,0,0)$, coefficients will be level 2 CMZVs. Also $$B = \frac{\splitfrac{\Gamma \left(d+\frac{1}{2}\right) \Gamma \left(e+\frac{1}{2}\right) \Gamma \left(-a+b+d+\frac{1}{2}\right) \Gamma \left(a-b+e+\frac{1}{2}\right)}{ \Gamma (-a+c+d) \Gamma (a-c+e+1) \Gamma \left(a-d+e+\frac{1}{2}\right) \Gamma \left(a-b-c-d+\frac{1}{2}\right)}}{\Gamma \left(\frac{1}{2}-b\right) \Gamma (1-c) \Gamma (a+e+1) \Gamma (b+e+1) \Gamma \left(c+e+\frac{1}{2}\right) \Gamma (d+e+1) \Gamma (a-b-c+1)}$$
	and $$C(n) = -\frac{(1-c)_n (a+e+1)_n (b+e+1)_n (d+e+1)_{2 n} (a-b-c+1)_n \times P}{\splitfrac{4 (n-c) (a+e+n) (b+e+n) (d+e+2 n-1) (d+e+2 n)}{ \left(d+\frac{1}{2}\right)_n \left(e+\frac{1}{2}\right)_n (a-b-c+n) \left(-a+b+d+\frac{1}{2}\right)_n \left(a-b+e+\frac{1}{2}\right)_n (a-c+e+1)_{2 n}}}$$
here $P\in \mathbb{Z}[a,b,c,d,e,n]$ whose full-form can be found in Mathematica notebook attached. \par 
Now comes third difficulty: when $(a,b,c,d,e)=(0,0,0,0,0)$, the series $\sum_{n\geq 1} C(n)$ diverges at rate $O(n^{-1})$: one can check
$$C(n) = \underbrace{\frac{-2^{-a+c+d} \Gamma \left(d+\frac{1}{2}\right) \Gamma \left(e+\frac{1}{2}\right)  \Gamma \left(-a+b+d+\frac{1}{2}\right) \Gamma \left(a-b+e+\frac{1}{2}\right)  \Gamma (a-c+e+1)}{\Gamma (1-c) \Gamma (a+e+1) \Gamma (b+e+1) \Gamma (d+e+1) \Gamma (a-b-c+1)}}_{:= D} n^{a-c-d-1} + O(n^{a-c-d-2}).$$
A way to remediate this is to write $$\sum_{n\geq 1}C(n) = \sum_{n\geq 1} \left(C(n) - D\times \Gamma(1+a-c-d)\frac{(1+a-c-d)_{n}}{(a-c-d+n)(1)_{n}}\right) + D\sum_{n\geq 1} \Gamma(1+a-c-d)\frac{(1+a-c-d)_{n}}{(a-c-d+n)(1)_{n}},$$
here the last sum is $-D\times\Gamma(a-c-d)$ if convergent, i.e. $\Re(a-c-d)<0$; while first sum is analytic under $\Re(a-c-d)<-1$. Thus we see that both sides of $$\sum_{k\geq 0}A(k) = \sum_{n\geq 1}\left(C(n) - \Gamma(1+a-c-d)\frac{(1+a-c-d)_{n}}{(a-c-d+n)(1)_{n}} \right) + B - D\Gamma(a-c-d)$$
are now analytic at $(a,b,c,d,e) = (0,0,0,0,0)$. Now we can finally compare coefficient of this equality on both sides, LHS will be CMZVs of level 2 (with weight = degree of monomial $+ 3$), the hypergeometric part of $\sum_{n\geq 1} C(n)$ is $1^n \frac{(1)_n^4}{(1/2)_n^4}$. 

For example, coefficient of constant term gives
$$\sum_{n\geq 1} \left( \frac{(1)_n^4}{(1/2)_n^4} \frac{1-4n}{4n^4} + \frac{\pi^2}{n}\right) = 14\zeta(3) -3\pi^2 \log 2;$$
coefficient of $a^1b^0c^0d^0e^0$ gives
$$\sum_{n\geq 1} \left( \frac{(1)_n^4}{(1/2)_n^4} \frac{(2-8 n) H_{n-1}+(4 n-1) H_{2 n-1}-1}{4 n^4}+\pi ^2 \frac{H_{n-1}-\log 2}{n}\right) = 64 \text{Li}_4\left(\frac{1}{2}\right)-\frac{53 \pi ^4}{90}+\frac{8 \log ^4(2)}{3}+\frac{29}{6} \pi ^2 \log ^2(2);$$
coefficient of $a^1b^1c^0d^0e^0$ gives
$$\sum_{n\geq 1} \left( \frac{(1)_n^4}{(1/2)_n^4}\frac{3 n (1-4 n) H_{n-1}^{(2)}+8 n (4 n-1) H_{2 n-1}^{(2)}+1}{4 n^5} -\frac{5\pi ^4}{6n}\right) = 7 \pi ^2 \zeta (3)-217 \zeta (5)+\frac{5}{2} \pi ^4 \log (2).$$
Taking a certain linear combination\footnote{specifically, $-[a^2] + [a b]/2 + [b^2] - [a d]/2 + [a e]/2$, here square bracket means coefficient of corresponding monomials.} of second degree monomials will cancel all occurence of harmonic numbers and produce $$\sum_{n\geq 1} \frac{(1)_n^4}{(1/2)_n^4}\frac{4 n-1}{4 n^5 (2 n-1)} = 93 \zeta (5)-7 \pi ^2 \zeta (3),$$
a conclusion we used in Example \ref{ex_zeta5}. 
\end{romexample}

\begin{remark}
One might ask whether a simpler method exists to prove the evaluation of  $\sum_{n\geq 1} \frac{(1)_n^4}{(1/2)_n^4}\frac{4 n-1}{4 n^5 (2 n-1)}$, this sum is actually a very-well-poised (and $2$-balanced) $_9F_8$: it equals $$12\times V(\frac{3}{2};\{1\}_6,\frac{1}{2}),$$
so a direct proof via Baliey's $4$-term transformation (\cite{bailey1935generalized}) formula is probable.
\end{remark}

\begin{romexample}\label{sun_100_prob_1}
Let $$F(n,k) = \frac{\splitfrac{4^n \Gamma \left(a+k+\frac{2}{3}\right) \Gamma \left(a+d+k+\frac{1}{3}\right) \Gamma (c+2 k-3 n+1)}{ \Gamma \left(2 a-2 b-c+n+\frac{4}{3}\right) (-2 a+b-d-2 k+n-1) \Gamma \left(2 a-b-c+d+2 n+\frac{1}{2}\right)}}{\splitfrac{\Gamma (b+n) \Gamma \left(a-b+k-n+\frac{5}{3}\right) \Gamma \left(-2 a+2 b+c-2 d-n+\frac{1}{3}\right)}{ \Gamma \left(2 a-2 b-c+d+n+\frac{1}{2}\right) \Gamma \left(a-b+d+k-n+\frac{4}{3}\right) \Gamma (4 a-2 b-c+2 d+2 k+n+2)}}.$$
We have $\mathcal{E}_F(x) = 4^2 | x-1| ^{2-2 x} | x| ^{2 x} | 2 x-3| ^{2 x-3} | 2 x+1| ^{-2 x-1}$, it does not satisfies the condition of Proposition \ref{vanishing_prop}, 
\begin{figure}[h]
	\centering
	\includegraphics[width=0.6\textwidth]{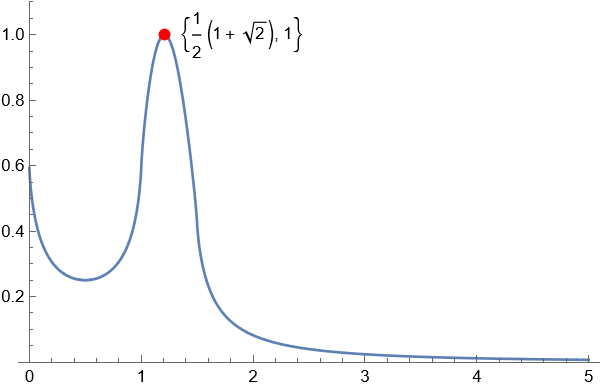}
	\caption{Plot of $\mathcal{E}_F(x)$ for Example \ref{sun_100_prob_1}, maximum is attained at $x=(1+\sqrt{2})/2$}
	\label{plot3}
\end{figure}
we will find $\lim_{n\to\infty}\sum_{k\geq 0} F(n,k)$ explicitly in lemma below. Proposition \ref{WZ_prop} becomes (with $g(n)=0$):
\begin{multline*}
b c (-6 a+6 b+3 c-1)  \sum_{k\geq 0} \frac{\left(a+\frac{2}{3}\right)_k (c+1)_{2 k} (2 a-b+d+2 k+1) (a+d+\frac{1}{3})_k}{(a-b+\frac{2}{3})_{k+1} (4 a-2 b-c+2 d+2 k+1) (a-b+d+\frac{1}{3})_{k+1} (4 a-2 b-c+2 d+1)_{2 k}} \\ =\frac{3\times 4^{-b} \Gamma (b+1) \Gamma (1-c) \Gamma \left(a-b+\frac{2}{3}\right) \Gamma \left(a-b+d+\frac{1}{3}\right) \Gamma \left(2 a-2 b-c+d+\frac{1}{2}\right) \Gamma (4 a-2 b-c+2 d+1)}{\Gamma \left(a+\frac{2}{3}\right) \Gamma \left(a+d+\frac{1}{3}\right) \Gamma \left(2 a-2 b-c+\frac{1}{3}\right) \Gamma \left(2 a-b-c+d+\frac{1}{2}\right) \Gamma \left(2 a-2 b-c+2 d+\frac{2}{3}\right)} \\ + \sum_{n\geq 0}\frac{\splitfrac{4^n (6 a-6 b-3 c+3 n+1)  P  (-a+b+\frac{1}{3})_n (2 a-2 b-c+\frac{1}{3})_n (-a+b-d+\frac{2}{3})_n}{ (2 a-b-c+d+\frac{1}{2})_{2 n} (2 a-2 b-c+2 d+\frac{2}{3})_n}}{9 (b+1)_n (1-c)_{3 n+2} (-4 a+4 b+2 c-2 d-2 n-1) (2 a-2 b-c+d+\frac{1}{2})_n (4 a-2 b-c+2 d+1)_n},
\end{multline*}
here $P\in \mathbb{Z}[a,b,c,n]$ whose full-form can be found in the Mathematica notebook attached. Coefficient of $a,b,c,d$ on LHS are level $3$ CMZVs. Letting $a=b=c=d=0$ we obtain $$\sum_{n\geq 0} \underbrace{(\frac{16}{27})^n \frac{(1/4)_n (1/3)_n (2/3)_n (3/4)_n}{(1/2)_n (1)_n^3}}_{a_n} \frac{66n^2+37n+4}{2n+1} = \frac{27\sqrt{3}}{2\pi}.$$
This was conjectured by Z.W. Sun in \cite{sun2023new}, where following extensions were also proposed, 
$$\sum_{n\geq 0} a_n \frac{(66 n^2+37 n+4) (-3 H_n-H_{2 n}+2 H_{4 n})+60 n+26}{2 n+1} = \frac{81\sqrt{3}\log 3}{2\pi},$$
$$\sum_{n\geq 0} a_n \frac{(66 n^2+37 n+4) (3 H_n-2 H_{2 n}+3 H_{3 n})-\frac{108 n^2+92 n+23}{2 n+1}}{2 n+1} = 0,$$
they are proved by comparing coefficients. 

\begin{lemma}
Let $F(n.k)$ be in above example, $a,b,c,d$ in a neighbourhood of origin, then
$$\lim_{n\to\infty}\sum_{k\geq 0} F(n,k) = 4^{-b} \csc (\pi  c) \sin (\pi (-2 a+2 b+c-2 d+\frac{1}{3})).$$
\end{lemma}
\begin{proof}[Proof sketch]
This can also be shown by using same technique as in Lemma \ref{lemma_hypertrans}, we introduce here yet another technique. Let $\alpha = (1+\sqrt{2})/2$, only $$\sum_{k\in (\alpha - 1/2, \alpha + 1/2)n} F(n,k)$$ contributes to the limit. Stirling's formula applies uniformly for $k$ in this range, let $1<\theta<3/2$, we have
\begin{multline*}F(n,\theta n) = \csc(\pi(1+c)) \sin(\pi(1/3-2a+2b+c-2d)) \times \exp\left[n f(\theta)\right] \\ \times \left(\frac{(1-2 \theta ) (3-2 \theta )^{c+\frac{1}{2}} \theta ^{2 a+d} 2^{2 a-b-c+d-\frac{1}{2}} (\theta -1)^{-2 a+2 b-d-2} (2 \theta +1)^{-4 a+2 b+c-2 d-\frac{3}{2}}}{\sqrt{n\pi }}\right) [1+ O(1/n)],\end{multline*}
here  $f(\alpha) = -2 \theta  \log (\theta -1)+2 \log (\theta -1)+2 \theta  \log (\theta )-2 \theta  \log (2 \theta +1)-\log (2 \theta +1)+2 \theta  \log (3-2 \theta )-3 \log (3-2 \theta )+4 \log (2)$. Letting $\theta = (\alpha n + r)/n$, then $$nf(\theta) =  -8r^2/n + O(r^3/n^2).$$ We then see it's allowed to further restrict the summation in $k$ to $k\in (\alpha n - n^{7/12}, \alpha n + n^{7/12})$,\footnote{actually any $n^\gamma$ with $2\gamma-1 >0$ and $3\gamma - 2 < 0$ works} expanding again at $n=\infty$ gives 
$$F(n,k) = \csc(\pi c)  \sin (\pi (-2 a+2 b+c-2 d+\frac{1}{3})) 4^{-b} \frac{2\sqrt{2}}{\sqrt{n\pi}} e^{-8r^2/n} \left[1 + O(n^{-1/4}) \right]\qquad k = \alpha n + r,$$
uniformly in this range of $k$, so we only need to find asymptotic of $$\sum_{r\in (-n^{7/12},n^{7/12})} e^{-8r^2/n}$$ as $n\to \infty$, this is easily done by first completing the sum to all of $r\in \mathbb{Z}$, which does not alter the leading asymptotic, then invoke Poisson summation formula. 
\end{proof}
\end{romexample}

\begin{romexample}\label{sun_100_prob_2}We give a final example in which all four terms in Proposition \ref{WZ_prop} are non-zero. This is motivated by the following conjecture of Z.W. Sun\footnote{\url{https://mathoverflow.net/questions/456528}}:
$$\sum_{n\geq 1} \left(\frac{3^4}{2^8}\right)^n \frac{(1/3)_n (2/3)_n (1)_n}{(1/4)_n (1/2)_n (3/4)_n} \frac{35n^2-29n+6}{n(3n-1)(3n-2)} = \sqrt{3}\pi.$$

Let $$F(n,k) = \frac{3^{4 n-3 k} \Gamma (a+3 c-2 n) \Gamma (-a+3 k-4 n+1)}{\Gamma (c+k-2 n+1) \Gamma (-a+3 b-3 c-n-1) \Gamma (2 a-3 b+6 c-n) \Gamma \left(-b+c+k-n+\frac{5}{3}\right) \Gamma \left(-a+b-2 c+k-n+\frac{4}{3}\right)}.$$
Now $\mathcal{E}_F(x) = \frac{1}{4} 3^{4-3 x} | x-2| ^{2-x} | x-1| ^{2-2 x} | 3 x-4| ^{3 x-4}$, 
\begin{figure}[h]
	\centering
	\includegraphics[width=0.6\textwidth]{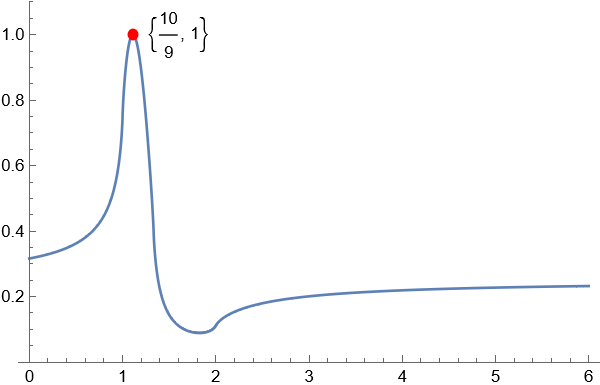}
	\caption{Plot of $\mathcal{E}_F(x)$ for Example \ref{sun_100_prob_2}, maximum is attained at $x=10/9$}
	\label{plot}
\end{figure}
it can shown, exactly as in previous lemma, that $$\lim_{n\to\infty} \sum_{k\geq 0} F(n,k) = -\pi^{-1}3^{2-a} \csc (\pi  a) \sin (\pi  c) \csc (\pi  (a+3 c)) \sin (\pi  (a-3 b+3 c)) \sin (\pi  (2 a-3 b+6 c)).$$
And
$$\sum_{k\geq 0} F(0,k) + \sum_{n\geq 0} g(n) = \lim_{n\to\infty}\sum_{k\geq 0} F(n,k) +  \sum_{n\geq 1} G(n-1,0) $$ becomes
\begin{multline*}
	\sum_{k\geq 0} \frac{3^{-3 k} (1-a)_{3 k}}{(c+1)_k (-b+c+\frac{2}{3})_{k+1} (-a+b-2 c+\frac{1}{3})_{k+1}}\\ - \frac{3^{\frac{5}{2}-a}\Gamma (c+1) \Gamma (-b+c+\frac{2}{3})\Gamma(-a+b-2 c+\frac{1}{3})}{2 \pi (a-3 b+3 c+1)  \Gamma (1-a)}\sum_{n\geq 0} \frac{  (-2 a+3 b-6 c+1)_n (a-3 b+3 c+1)_{n+1}}{ (-a-3 c+1)_{2 n+2}} \\ = \frac{3^{2-a} \Gamma (a) \Gamma (-a-3 c+1) \Gamma \left(-b+c+\frac{2}{3}\right) \Gamma \left(-a+b-2 c+\frac{1}{3}\right)}{\Gamma (-c) \Gamma (-2 a+3 b-6 c+1) \Gamma (a-3 b+3 c+2)} \\ + \frac{c}{a(a-3 b+3 c+1)}\sum_{n\geq 1} \frac{3^{4 n}  \times P (1-c)_{2 n-2} (b-c+\frac{1}{3})_{n-1} (-2 a+3 b-6 c+1)_{n-1} (a-b+2 c+\frac{2}{3})_n (a-3 b+3 c+1)_n}{(a+1)_{4 n-1} (3 a-3 b+6 c+3 n-1) (-a-3 c+1)_{2 n}},
\end{multline*}
where $P = a^3+6 a^2 c-27 a b c+54 a b n-27 a b+36 a c^2-69 a c n+42 a c+42 a n^2-57 a n+17 a+27 b^2 c-54 b^2 n+27 b^2-81 b c^2+162 b c n-90 b c+18 b n-9 b+54 c^3-99 c^2 n+72 c^2-21 c n^2-30 c n+18 c+70 n^3-93 n^2+41 n-6$.
The case $a=b=c=0$ reduces to the conjecture.

\end{romexample}
~\\[0.05in]
\textbf{Acknowledgements} The author would like to express his sincere gratitude to W. Zudilin, D. Zagier and Z.W. Sun for valuable discussions.
	
\bibliographystyle{plain} 
\bibliography{ref.bib} 

\begin{thebibliography}{10}

\bibitem{au2022multiple}
Kam~Cheong Au.
\newblock Multiple zeta values, {WZ}-pairs and infinite sums computations.
\newblock {\em arXiv preprint arXiv:2212.02986}, 2022.

\bibitem{au2024wilfQ}
Kam~Cheong Au.
\newblock Wilf-{Z}eilberger seeds: $ q $-analogues.
\newblock {\em arXiv preprint arXiv:2403.04555}, 2024.

\bibitem{bailey1935generalized}
Wilfrid~Norman Bailey.
\newblock Generalized hypergeometric series.
\newblock 1935.

\bibitem{borwein1987pi}
Jonathan~M Borwein and Peter~B Borwein.
\newblock {\em Pi and the AGM: a study in the analytic number theory and
  computational complexity}.
\newblock Wiley-Interscience, 1987.

\bibitem{broadhurst2024five}
David Broadhurst.
\newblock Five families of rapidly convergent evaluations of zeta values.
\newblock {\em arXiv preprint arXiv:2401.08997}, 2024.

\bibitem{chu2011dougall}
Wenchang Chu.
\newblock Dougall’s bilateral {$_2H_2$}-series and {R}amanujan-like
  {$\pi$}-formulae.
\newblock {\em Mathematics of computation}, 80(276):2223--2251, 2011.

\bibitem{chu2014accelerating}
Wenchang Chu and Wenlong Zhang.
\newblock Accelerating {D}ougalls {$_5F_4$}-sum and infinite series involving
  $\pi$.
\newblock {\em Mathematics of Computation}, 83(285):475--512, 2014.

\bibitem{cohen2021rational}
Henri Cohen and Jes{\'u}s Guillera.
\newblock Rational hypergeometric {R}amanujan identities for $1/\pi^c$: Survey
  and generalizations.
\newblock {\em arXiv preprint arXiv:2101.12592}, 2021.

\bibitem{de1981asymptotic}
Nicolaas~Govert De~Bruijn.
\newblock {\em Asymptotic methods in analysis}, volume~4.
\newblock Courier Corporation, 1981.

\bibitem{dembele2022special}
Lassina Demb{\'e}l{\'e}, Alexei Panchishkin, John Voight, and Wadim Zudilin.
\newblock Special hypergeometric motives and their {L}-functions: Asai
  recognition.
\newblock {\em Experimental Mathematics}, 31(4):1278--1290, 2022.

\bibitem{gessel1995finding}
Ira~M Gessel.
\newblock Finding identities with the {WZ} method.
\newblock {\em Journal of symbolic computation}, 20(5-6):537--566, 1995.

\bibitem{guillera2003new}
Jes{\'u}s Guillera.
\newblock About a new kind of {R}amanujan-type series.
\newblock {\em Experimental Mathematics}, 12(4):507--510, 2003.

\bibitem{guillera2008hypergeometric}
Jes{\'u}s Guillera.
\newblock Hypergeometric identities for 10 extended {R}amanujan-type series.
\newblock {\em The Ramanujan Journal}, 15(2):219--234, 2008.

\bibitem{guillera2011new}
Jes{\'u}s Guillera.
\newblock A new ramanujan-like series for 1/$\pi$ 2.
\newblock {\em The Ramanujan Journal}, 3(26):369--374, 2011.

\bibitem{guillera2016bilateral}
Jes{\'u}s Guillera.
\newblock Bilateral sums related to {R}amanujan-like series.
\newblock {\em arXiv preprint arXiv:1610.04839}, 2016.

\bibitem{guillera2018dougall}
Jes{\'u}s Guillera.
\newblock Dougall’s {$_5F_4$} sum and the {WZ} algorithm.
\newblock {\em The Ramanujan Journal}, 46:667--675, 2018.

\bibitem{hou2023taylor}
Qing-Hu Hou and Zhi-Wei Sun.
\newblock Taylor coefficients and series involving harmonic numbers.
\newblock {\em arXiv preprint arXiv:2310.03699}, 2023.

\bibitem{mohammed2005infinite}
Mohamud Mohammed.
\newblock Infinite families of accelerated series for some classical constants
  by the {M}arkov-{WZ} method.
\newblock {\em Discrete Mathematics \& Theoretical Computer Science}, 7, 2005.

\bibitem{paule1995mathematica}
Peter Paule and Markus Schorn.
\newblock A mathematica version of zeilberger's algorithm for proving binomial
  coefficient identities.
\newblock {\em Journal of symbolic computation}, 20(5-6):673--698, 1995.

\bibitem{pilehrood2008generating}
Kh~Hessami Pilehrood and T~Hessami Pilehrood.
\newblock Generating function identities for $\zeta(2n+ 2)$,$\zeta(2n+ 3)$ via
  the {WZ} method.
\newblock {\em The Electronic Journal of Combinatorics}, pages R35--R35, 2008.

\bibitem{pilehrood2008simultaneous}
Tatiana~Hessami Pilehrood and Khodabakhsh~Hessami Pilehrood.
\newblock Simultaneous generation for zeta values by the {M}arkov-{WZ} method.
\newblock {\em Discrete Mathematics \& Theoretical Computer Science}, 10, 2008.

\bibitem{pilehrood2010series}
Tatiana~Hessami Pilehrood and Khodabakhsh~Hessami Pilehrood.
\newblock Series acceleration formulas for beta values.
\newblock {\em Discrete Mathematics \& Theoretical Computer Science}, 12, 2010.

\bibitem{shimura1971introduction}
Gor{\=o} Shimura.
\newblock {\em Introduction to the arithmetic theory of automorphic functions},
  volume~1.
\newblock Princeton university press, 1971.

\bibitem{sun2021book}
Zhi-Wei Sun.
\newblock New conjectures in number theory and combinatorics.
\newblock {\em Harbin Institute of Technology Press, Harbin}, 2021.

\bibitem{sun2022conjectures}
Zhi-Wei Sun.
\newblock Series with summands involving higher order harmonic numbers.
\newblock {\em arXiv preprint arXiv:2210.07238}, 2022.

\bibitem{sun2023new}
Zhi-Wei Sun.
\newblock New series involving binomial coefficients.
\newblock {\em arXiv preprint arXiv:2307.03086}, 2023.

\bibitem{wei2022conjectural}
Chuanan Wei.
\newblock On a conjectural series for {$\pi$} and its {$q $}-analogue.
\newblock {\em arXiv preprint arXiv:2211.11484}, 2022.

\bibitem{wei2023conjectural}
Chuanan Wei.
\newblock On a conjectural series of sun for the mathematical constant
  {$\beta(4)$}.
\newblock {\em arXiv preprint arXiv:2303.05402}, 2023.

\bibitem{wei2023some2}
Chuanan Wei.
\newblock On some conjectural series containing binomial coefficients and
  harmonic numbers.
\newblock {\em arXiv preprint arXiv:2306.02641}, 2023.

\bibitem{wei2023some}
Chuanan Wei.
\newblock Some fast convergent series for the mathematical constants $\zeta(4)$
  and $\zeta(5)$.
\newblock {\em arXiv preprint arXiv:2303.07887}, 2023.

\bibitem{wei2023some3}
Chuanan Wei and Ce~Xu.
\newblock On some conjectural series containing harmonic numbers of 3-order.
\newblock {\em arXiv preprint arXiv:2308.06440}, 2023.

\bibitem{wei2023two}
Chuanan Wei and Ce~Xu.
\newblock On two conjectural series involving riemann zeta function.
\newblock {\em arXiv preprint arXiv:2310.04642}, 2023.

\bibitem{zhang2015common}
Wenlong Zhang.
\newblock Common extension of the {W}atson and {W}hipple sums and
  {R}amanujan-like $\pi$-formulae.
\newblock {\em Integral Transforms and Special Functions}, 26(8):600--618,
  2015.

\bibitem{zhao2016multiple}
Jianqiang Zhao.
\newblock {\em Multiple zeta functions, multiple polylogarithms and their
  special values}, volume~12, chapter~13.
\newblock World Scientific, 2016.

\bibitem{zudilin2007quadratic}
Wadim Zudilin.
\newblock Quadratic transformations and {G}uillera’s formulas for 1/$\pi^2$.
\newblock {\em Mathematical Notes}, 81, 2007.

\bibitem{zudilin2011arithmetic}
Wadim Zudilin.
\newblock Arithmetic hypergeometric series.
\newblock {\em Russian Mathematical Surveys}, 66(2):369, 2011.

\end{thebibliography}
	
\end{document}